\documentclass[UTF-8,reqno]{amsart}
\usepackage{enumerate}
\usepackage{mhequ}
\usepackage{amsthm,amsmath,amssymb,url,color, booktabs,nccmath}
\usepackage[left=3.2cm,right=3.2cm,top=4cm,bottom=4cm]{geometry}
\usepackage{mathrsfs}
\usepackage{enumitem,dsfont}
\usepackage{subfigure}
\usepackage[graphicx]{realboxes}

\usepackage{color}
\usepackage[colorlinks=true]{hyperref}
\hypersetup{
    linkcolor=blue,          
    citecolor=red,        
    filecolor=blue,      
    urlcolor=cyan
}

\definecolor{darkergreen}{rgb}{0.0, 0.5, 0.0}

\numberwithin{equation}{section}
\def\theequation{\arabic{section}.\arabic{equation}}
\newcommand{\be}{\begin{eqnarray}}
\newcommand{\ee}{\end{eqnarray}}
\newcommand{\ce}{\begin{eqnarray*}}
\newcommand{\de}{\end{eqnarray*}}
\newtheorem{theorem}{Theorem}[section]
\newtheorem{lemma}[theorem]{Lemma}
\newtheorem{proposition}[theorem]{Proposition}
\newtheorem{Examples}[theorem]{Example}
\newtheorem{corollary}[theorem]{Corollary}

\newtheorem{definition}[theorem]{Definition}
\theoremstyle{definition}
\newtheorem{remark}[theorem]{Remark}

\DeclareMathOperator{\supp}{supp}

\def\u{\mathbf{u}}

\def\[{{\Big[}}
\def\]{{\Big]}}
\def\<{{\langle}}
\def\>{{\rangle}}
\def\({{\Big(}}
\def\){{\Big)}}

\def\bx{{\mathbf{x}}}
\def\tr{\mathrm {tr}}

\def\dif{{\mathord{{\rm d}}}}

\def\={&\!\!=\!\!&}

\def\mN{{\mathbb N}}

\def\bP{{\mathbf P}}

\def\1{{\mathbf{1}}}

\def\geq{\geqslant}
\def\leq{\leqslant}

\def\div{\mathord{{\rm div}}}

\def\u{\mathbf{u}}

\def\[{{\Big[}}
\def\]{{\Big]}}
\def\<{{\langle}}
\def\>{{\rangle}}
\def\({{\Big(}}
\def\){{\Big)}}

\def\bx{{\mathbf{x}}}
\def\tr{\mathrm {tr}}

\def\dif{{\mathord{{\rm d}}}}

\def\={&\!\!=\!\!&}
\def\bt{\begin{theorem}}
\def\et{\end{theorem}}
\def\bl{\begin{lemma}}
\def\el{\end{lemma}}
\def\br{\begin{remark}}
\def\er{\end{remark}}
\def\bx{\begin{Examples}}
\def\ex{\end{Examples}}
\def\bd{\begin{definition}}
\def\ed{\end{definition}}
\def\bp{\begin{proposition}}
\def\ep{\end{proposition}}
\def\bc{\begin{corollary}}
\def\ec{\end{corollary}}

\def\geq{\geqslant}
\def\leq{\leqslant}

\def\div{\mathord{{\rm div}}}

\def\Id{\textrm{Id}}

\def\bP{{\mathbf P}}

 \def\R{\mathbb R}
 \def\R{\mathbb R}    
\def\N{\mathbb N}  
   
\def\<{\langle} \def\>{\rangle}

\allowdisplaybreaks

\begin{document}

\title[PROBABILISTICALLY STRONG SOLUTIONS TO STOCHASTIC POWER-LAW EQUATION]{Global-in-time PROBABILISTICALLY STRONG SOLUTIONS TO STOCHASTIC POWER-LAW EQUATIONS: EXISTENCE AND NON-UNIQUENESS}

\author{Huaxiang L\"u}
\address[H. L\"u]{Academy of Mathematics and Systems Science,
Chinese Academy of Sciences, Beijing 100190, China}
\email{lvhuaxiang22@mails.ucas.ac.cn }

\author{Xiangchan Zhu}
\address[X. Zhu]{ Academy of Mathematics and Systems Science,
Chinese Academy of Sciences, Beijing 100190, China; Fakult\"at f\"ur Mathematik, Universit\"at Bielefeld, D-33501 Bielefeld, Germany}
\email{zhuxiangchan@126.com}

\thanks{
Research  supported   by National Key R\&D Program of China (No. 2020YFA0712700) and the NSFC (No.  12090014, 12288201) and
  the support by key Lab of Random Complex Structures and Data Science,
 Youth Innovation Promotion Association (2020003), Chinese Academy of Science. The financial support by the DFG through the CRC 1283 ``Taming uncertainty and profiting
 from randomness and low regularity in analysis, stochastics and their applications'' is greatly acknowledged.
}

\begin{abstract}
We are concerned with the power-law fluids
driven by an additive stochastic forcing in dimension $d\geq3$. For the power index $r\in(1,\frac{3d+2}{d+2})$, we establish existence of infinitely many global-in-time probabilistically strong
and analytically weak solutions in $L^p_{loc}([0,\infty);L^2)\cap C([0,\infty);W^{1,\max\{1,r-1\}}),p\geq1$ for every divergence free initial condition in $L^2\cap W^{1,\max\{1,r-1\}}$. This result in
particular implies non-uniqueness in law.  Our result is sharp in the three dimensional case in the sense that the solution is unique if $r\geq \frac{3d+2}{d+2}$.
\end{abstract}


\date{\today}

\maketitle

\tableofcontents

\section{Introduction}

In this paper we are concerned with the stochastic power-law flows on $\mathbb{T}^d\ (d\geq3)$ driven by an additive noise. The equations read as
\begin{align}\label{eq:pl}
 \dif u+\div(u\otimes u)\dif t-\div\mathcal{A}(Du)\dif t+\nabla{\pi}\dif t&=\dif B,\notag\\
\div u&=0, \\
u(0)&=u_0,\notag
\end{align}
where $\pi$ is the associated pressure, $Dv=\frac12(\nabla v+\nabla^T v)$, and the
non-Newtonian tensor $\mathcal{A}$ is given by the following power law
\begin{align}
\mathcal{A}(Q)=(\nu_0+\nu_1|Q|)^{r-2}Q,\label{AQ1}
\end{align}
for some $\nu_0\geq 0,\nu_1\geq 0$ and $r\in(1, \infty)$. Here $B$ is a $GG^*$-Wiener process on some probability space $(\Omega, \mathcal{F},\mathbf{P})$ and $G $ is a Hilbert-Schmidt operator from $U$ to $L^2$ for some Hilbert space $U$. In addition we require certain regularity of the noise, namely we assume there exists a $\sigma>0$ such that $\mathrm{Tr}((-\Delta)^{2\sigma}GG^*)<\infty$ (see Remark \ref{jieshiz} below).

The power-law flows have been proposed independently by Norton \cite{NF} in metallurgy and by de Waele \cite{DW} and Ostwald \cite{O} in polymer chemistry. The case $r=2$ is the Navier-Stokes equations as $\div(Dv)=\frac{1}{2}\div(\nabla v+\nabla^Tv)=\frac{1}{2}(\Delta v+\nabla\div v)=\frac{1}{2}\Delta v$. When $r>2$, (\ref{eq:pl}) models the shear thickening fluids, whereas $r<2$ (\ref{eq:pl}) describes shear thinning fluids. Both cases have many applications in science and engineering, such as automobile engine oil, body armors and so on.

The mathematical discussion of power-law flows started from Lions and Ladyshenskaya \cite{L2,L3,L1,L4}.
 They proved that when $r\geq \frac{3d+2}{d+2}$ there exists a weak solution in the space
$L^r([0, T] ; W_{0,\div}^{1,r}(\mathbb{T}^d))\cap L^\infty([0, T] ;L^2(\mathbb{T}^d))$, where $W_{0,\div}^{1,r}(\mathbb{T}^d)$ is the closure of smooth functions on $\mathbb{T}^d$  and divergence-free in $W^{1,r}(\mathbb{T}^d)$. 
 Wolf \cite{W} improved this result to the case $r>\frac{2d+2}{d+2} $ using $L^\infty$-truncation. Then it was improved to $r>\frac{2d}{d+2}$ by Diening, R$\rm \mathring{u}\check{z}$i$\rm\check{c}$ka and Wolf \cite{LMJ} via the Lipschitz truncation method. Ko \cite{Ko22} consider the variable power-law index and  investigate the decay
properties of strong solutions based on the Fourier splitting method. Concerning  the uniqueness of solution, Ladyshenskaya \cite{L1} established the uniqueness provided $r\geq{1+\frac{d}{2}}$ or in case of smooth initial condition for $r\geq\frac{3d+2}{d+2}$.
M\'alek, Ne$\rm\check{c}$as, Rokyta and R$\rm \mathring{u}\check{z}$i$\rm\check{c}$ka  \cite{MNRR} proved the uniqueness for initial condition in $H^1$ provieded $r\geq\frac{3d+2}{d+2}$. More recently, by using the method of convex integration Burczak, Modena and Sz\'ekelyhidi Jr. \cite{BMS} obtained  the existence of multiple solutions for $r\in(1,\frac{3d+2}{d+2})$.

Unlike the deterministic counterpart, the stochastic equations possess additional structural features. We distinguish between probabilistically strong and probabilistically weak (martingale) solutions.
Probabilistically strong solutions are constructed on a given probability space and are adapted with respect to the given noise, while probabilistically weak solutions are typically obtained by the method of compactness where the noise as well as the probability space becomes part of the construction.
Using the method of compactness in the stochastic setting Terasawa and Yoshida \cite{YN,Yos12} gave the existence of probabilistically weak solutions in the case $\nu_0>0$ for
\begin{align*}
r\in
\begin{cases}
(\frac{3d}{d+2}\vee\frac{3d-4}{d},\infty),&2\leq d\leq8,\\
(\frac{23}{9},\frac{18}{7})\cup(\frac{19+\sqrt{793}}{18},\infty),&d=9,\\
(\frac{3d-8+\sqrt{9d^2+64}}{2d},\infty),&d\geq10.
\end{cases}
\end{align*}
 Breit \cite{BD} showed the existence of probabilistically  weak solutions provided $r>\frac{2d+2}{d+2}$ in a
bounded Lipschitz domain based on the $L^\infty$-truncation and a harmonic pressure decomposition
which are adapted to the stochastic setting. When $r\geq1+\frac{d}{2}$ Terasawa and Yoshida \cite{YN,Yos12} also obtained pathwise uniqueness which by the Yamada-Watanabe theorem implies the existence and uniqueness of probabilistically strong solution.

However, it is not clear that  probabilistically strong solutions exist when $r\in(1, 1+\frac{d}{2})$ which has more practical significance because many different fluids lie in the interval. In fact, it is necessary to take expectation to control the noise and obtain uniform estimates, which then leads to probabilistically weak solutions. When $r \in(\frac{2d+2}{d+2}, 1+\frac{d}{2})$ due to lack of uniqueness we cannot apply the Yamada-Watanabe theorem to obtain probabilistically strong solutions. If we analyze the equation $\omega$-wise, then the converging subsequence from compactness argument may depend on $\omega$ which destroys adaptedness. Moreover, when   $r\in(1,\frac{2d+2}{d+2}]$ there's even no result for the existence of probabilistically weak solution. The first goal of the present paper is to  establish the existence of global-in-time  probabilistically strong and analytically weak solutions when $r\in(1,\frac{3d+2}{d+2})$.

Another natural question of stochastic system is that whether the uniqueness in law holds when $r\in (1,1+\frac{d}{2})$. In the stochastic setting, there are explicit examples of stochastic differential equations where pathwise uniqueness does not hold but uniqueness in law is valid. Moreover, there is evidence that a suitable stochastic perturbation may provide a regularizing effect on deterministic ill-posed problems. A transport noise and linear multiplicative
noise prevent blow up of strong solutions have been obtained by Flandoli, Gubinelli and Priola \cite{FGP10} and
Flandoli and Luo \cite{FL19} and Glatt-Holtz and Vicol \cite{GHV14} and R\"ockner, Zhu and Zhu \cite{RZZ14}. In the present paper we prove non-uniqueness in law holds when $r\in(1,\frac{3d+2}{d+2})$. This result is sharp in the three dimensional case in the sense that there exists a unique global probabilistically strong solution when $r\geq\frac{3d+2}{d+2}=\frac{11}{5}$.

\subsection{Main result}

Our first result is the existence of global-in-time non-unique probabilistically strong and analytically weak solutions to (\ref{eq:pl}) for every given divergence free initial condition in $L^2\cap W^{1,\max\{1,r-1\}}$.

\bt\label{thm:4}
Let $r\in(1,\frac{3d+2}{d+2}), r^*= \max\{1,r-1\}$ and  $s\in[r^*,\frac{2d}{d+2})$. Assume $u_0\in L^2\cap W^{1,s}\ \mathbf{P}$-a.s. is a divergence free initial condition independent of the Wiener process B. There exist infinitely many probabilistically strong and analytically weak solutions to the power-law system (\ref{eq:pl}) on $[0, \infty)$. The solutions belong to $L^p_{loc}([0, \infty);L^2)\cap C([0, \infty); W^{1,s})\ \mathbf{P}$-a.s. for all $p\in[1, \infty).$
\et

Furthermore, our construction directly implies the following result. 
\bc\label{thm:10}
Let $r\in(1,\frac{3d+2}{d+2})$, non-uniqueness in law holds for the power-law system (\ref{eq:pl}) for every given
initial law supported on divergence free vector fields in $L^2\cap W^{1,\max\{1,r-1\}}$.
\ec
\br
 When $r=2$ the power-law flows are the Navier-Stokes equations. Hofmanov\'{a}, Zhu, and the second named author \cite{HZZ} proved that
 for the initial value $u_{0}\in L^2$ $\mathbf{P}$-a.s. there exist infinitely many global strong solutions in $L^p_{\rm loc}([0,\infty);L^2)\cap C((0,\infty);W^{\frac12,\frac{31}{30}})$ $\mathbf{P}$-a.s. for all $p\in[1,\infty)$.
Compared to this, the solutions obtained in Theorem \ref{thm:4} have better regularity provided the initial value have better regularity.
\er

We prove Theorem \ref{thm:4} by the convex integration method in the stochastic setting developed in \cite{HZZ}.  In particular, we could construct a convex integration solution before a stopping times, which has been achieved for deterministic counterpart, see \cite{BMS}. Compared to \cite{BMS},  we use the intermittent jets
from \cite{BCV18} which seems to be better suited for the stochastic setting, see Remark \ref{BMSno}. To this end we extend the intermittent jets to higher dimensions. Moreover, as the operator $\mathcal{A}$ might be degenerate there's no  smoothing property as the heat semigroup which is required in the convex integration scheme for the initial value in \cite{HZZ}. In the following to solve this problem we still introduced the Laplacian term in the linear equation and subtract the corresponding term in the nonlinear equation (see \eqref{linear} and \eqref{eq:v1} below for more details).  
 And then we could extend the convex integration solution beyond a stopping time by connecting it to another strong solution and finally obtain the global-in-time solution.

Our result is sharp in 3D case in the sence that if $r\geq\frac{3d+2}{d+2}=\frac{11}{5}$, then the solution is unique.
\bt\label{thm:d=3}
Let $d=3,\nu_0>0$, and $G$ satisfies $\|\nabla G\|_{L_2(U,L_\sigma^2)}<\infty$. Then for $r\geq \frac{11}{5}, u_0\in H^1$, there exists a unique global probabilistically strong solution to (\ref{eq:pl}) in $C([0,T];L_\sigma^2)\cap L^r([0,T];W^{1,r})$ to (\ref{eq:pl}).
\et

Now we summarize existence and uniqueness/non-uniqueness of the solutions to (\ref{eq:pl}) as follows.
\includegraphics[width=14cm]{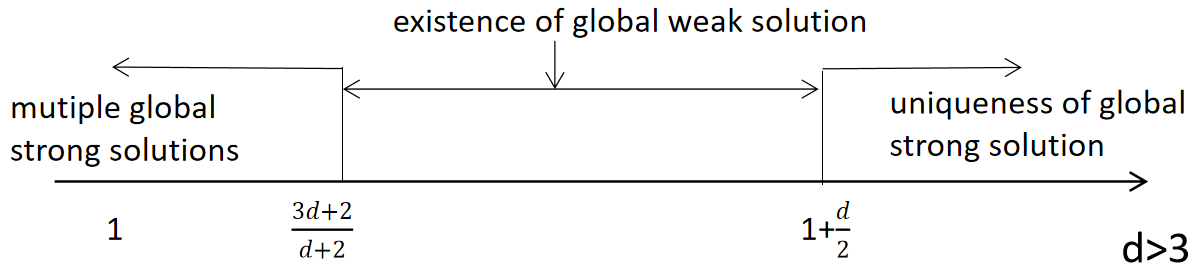}

\includegraphics[width=15cm]{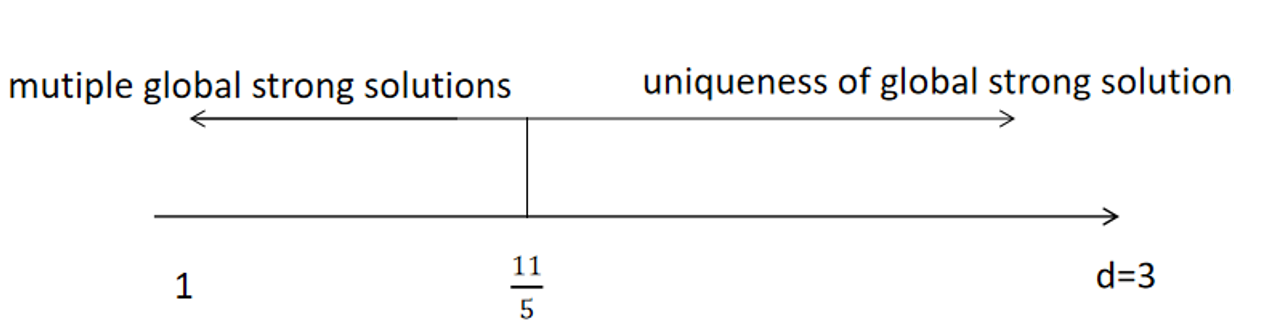}

\subsection{Further relevant literature}
The convex integration was introduced to fluid dynamics by De Lellis and Sz\'ekelyhidi Jr. \cite{DLS09, DLS10, DLS13}. This method has already led to a number of groundbreaking results concerning the incompressible Euler equations, culminating in the proof of Onsager’s conjecture by Isett \cite{Ise18} and by Buckmaster, De Lellis, Sz\'ekelyhidi Jr. and Vicol \cite{BDLSV19}. Also the question of well/ill-posedness of the three dimensional Navier-Stokes
equations has experienced an immense breakthrough: Buckmaster and Vicol \cite{BVb} established non-uniqueness of weak solutions with finite kinetic energy. More precisely, the authors
showed that for any prescribed smooth and non-negative function $e$ there is a weak solution whose
kinetic energy is given by $e$. Remarkably, Buckmaster, Colombo and Vicol \cite{BCV18} were even able
to connect two arbitrary strong solutions via a weak solution. Sharp non-uniqueness results for the Navier–Stokes equations
in dimension $d \geq2$ were obtained by Cheskidov and Luo \cite{AX}. The interested
reader is referred to the nice reviews \cite{BV,BV21} for further details and references. By a different method, a first non-uniqueness result for Leray solutions was established by Albritton, Bru\'{e} and Colombo \cite{ABC21} for the Navier-Stokes system with a force. 
An adaptation of convex integration to the stochastic setting has already appeared in a number of works, proving results of ill-posedness in various settings, see \cite{CFF19,CDZ22, BFH20c,
HZZ19, HZZ20, HZZ,HZZ22,Yam20a, Yam20b, Yam21a, Yam21b, Yam21c, RS21}.  \\

\noindent{\bf Organization of the paper.}
In Section~\ref{s:not} we collect the basic notations used throughout the~paper. In Section \ref{gij} we introduce the generalized intermittent jets, which play a crucial role in the construction of perturbation $w_{q+1}$.  Section \ref{cogpss} is devoted to our main convex integration result, Theorem \ref{thm:4} and Corollary \ref{thm:10}.   The proof of Theorem \ref{thm:d=3} is presented in Section \ref{tctr}. In Appendix \ref{s:appA} we give estimates of amplitude functions used in the convex integration construction. In Appendix \ref{tec} we give some technical tools used in the proof.

\section{Notations}
\label{s:not}

  Throughout the paper, we employ the notation $a\lesssim b$ if there exists a constant $c>0$ such that $a\leq cb$, and we write $a\simeq b$ if $a\lesssim b$ and $b\lesssim a$. $\mN_{0}:=\mN\cup \{0\}$. Given a Banach space $E$ with a norm $\|\cdot\|_E$ and $T>0$, we write $C_TE=C([0,T];E)$ for the space of continuous functions from $[0,T]$ to $E$, equipped with the supremum norm $\|f\|_{C_TE}=\sup_{t\in[0,T]}\|f(t)\|_{E}$. We also use $CE$ or $C([0,\infty);E)$ to denote the space of continuous functions from $[0,\infty)$ to $E$. For $\alpha\in(0,1)$ we  define $C^\alpha_TE$ as the space of $\alpha$-H\"{o}lder continuous functions from $[0,T]$ to $E$, endowed with the norm $\|f\|_{C^\alpha_TE}=\sup_{s,t\in[0,T],s\neq t}\frac{\|f(s)-f(t)\|_E}{|t-s|^\alpha}+\sup_{t\in[0,T]}\|f(t)\|_{E}.$ Here we use $C_T^\alpha$ to denote the case when $E=\mathbb{R}$. We also use $C_\mathrm{loc}^\alpha E$ to denote the space of functions from $[0,\infty)$ to $E$ satisfying $f|_{[0,T]}\in C_T^\alpha E$ for all $T>0$. For $p\in [1,\infty]$ we write $L^p_TE=L^p([0,T];E)$ for the space of $L^p$-integrable functions from $[0,T]$ to $E$, equipped with the usual $L^p$-norm. We also use $L^p_{\mathrm{loc}}([0,\infty);E)$ to denote the space of functions $f$ from $[0,\infty)$ to $E$ satisfying $f|_{[0,T]}\in L^p_T E$ for all $ T>0$.
    We use $L^p$ to denote the set of  standard $L^p$-integrable functions from $\mathbb{T}^d$ to $\mathbb{R}^d$. For $s>0$, $p>1$ we set $W^{s,p}:=\{f\in L^p; \|(I-\Delta)^{\frac{s}{2}}f\|_{L^p}<\infty\}$ with the norm  $\|f\|_{W^{s,p}}=\|(I-\Delta)^{\frac{s}{2}}f\|_{L^p}$. Set $L^{2}_{\sigma}=\{f\in L^2; \int_{\mathbb{T}^{d}} f\,\dif x=0,\div f=0\}$. For $s>0$, we define $H^s:=W^{s,2}\cap L^2_\sigma$.
For $T>0$ and a domain $D\subset\R^{+}$ we denote by  $C^{N}_{T,x}$ and $C^{N}_{D,x}$, respectively, the space of $C^{N}$-functions on $[0,T]\times\mathbb{T}^{d}$ and on $D\times\mathbb{T}^{d}$, respectively,  $N\in\N_{0}$. The spaces are equipped with the norms
$$
\|f\|_{C^N_{T,x}}=\sum_{\substack{0\leq n+|\alpha|\leq N\\ n\in\N_{0},\alpha\in\N^{d}_{0} }}\|\partial_t^n D^\alpha f\|_{L^\infty_T L^\infty},\qquad \|f\|_{C^N_{D,x}}=\sum_{\substack{0\leq n+|\alpha|\leq N\\ n\in\N_{0},\alpha\in\N^{d}_{0} }}\sup_{t\in D}\|\partial_t^n D^\alpha f\|_{ L^\infty}.
$$
We define the projection onto null-mean functions
is $\mathbb{P}_{\neq0} f := f -(2\pi)^{-d}\int_{\mathbb{T}^d} f\dif x$. For a matrix $T$, we denote its traceless part by $\mathring{T} := T- \frac{1}{d} \tr (T )\Id$. By $\mathcal{S}^{d\times d}$ we denote the space of symmetric matrix and by $\mathcal{S}_0^{d\times d}$ the space of symmetric trace-free matrix.

Regarding the driving noise, we assume that $B$ is a $GG^*$-Wiener process on some probability space $(\Omega, \mathcal{F}, \mathbf{P})$ and $G$ is a Hilbert--Schmidt operator from $U$ to $L_{\sigma}^2(\mathbb{T}^{d})$ for some Hilbert space $U$.
In addition we assume there exists a $\sigma>0$ such that $\mathrm{Tr}((-\Delta)^{2\sigma}GG^*)<\infty$.

\section{Generalized intermittent jets}\label{gij}
In this section we generalize the intermittent jets introduced in \cite{BV, BCV18} to higher dimensions. We point out that the construction is entirely deterministic, that is, none of the functions below depends on $\omega$. In \cite{BV, BCV18}   the authors used the $curl$ operator in the incompressibility corrector. However, there's no $curl$ operator in higher dimensions. Instead we introduced the corrector term inspired by \cite{AX}.

Let us begin with the following geometric lemma which can be found in \cite[Lemma 4.2]{AX}.
\bl\label{lem:3}
Denote by $\bar{B}_{1/2}(\mathrm{Id})$ the closed ball of radius $1/2$ around the identity matrix $\mathrm{Id}$, in the
space of $d\times d$ symmetric matrices. There exists $\Lambda\in \mathbb{S}^{d-1}\cap \mathbb{Q}^d$ such that for each $\xi\in\Lambda$ there exists
a $C^\infty$-function $\gamma_\xi$: $\bar{B}_{1/2}(\mathrm{Id})\to \mathbb{R}$ such that
$$R=\sum_{\xi\in\Lambda}\gamma_\xi^2(R)(\xi\otimes\xi)$$
for every symmetric matrix satisfying $|R-\mathrm{Id}|\leq 1/2$.
For $C_\Lambda=8|\Lambda|(1+(2\pi)^d)^{1/2}$, where $|\Lambda|$ is the cardinality of $\Lambda$, we define the constant
$$M^*=C_\Lambda\sum_{\xi\in\Lambda}(\|\gamma_\xi\|_{C^0}+\sum_{|j|\leq N^*}\|D^j\gamma_\xi\|_{C^0}).$$
\el
Here $N^*=[\frac{dN}2]+2$ is introduced in (\ref{N*}).

For parameters $\lambda,r_\perp, r_\parallel > 0 $, we assume
$$\lambda^{-1}\ll r_{\perp}\ll r_{\parallel}\ll 1,\ \ \lambda r_{\perp}\in\mathbb{N}.$$
\bl\label{lem:4}
$($\cite[Lemma 3]{BMS}$)$
Let $d\geq3$. Then there exist $\alpha_\xi$ and $\rho>0$ such that
$$(B_\rho(\alpha_\xi)+\{s\xi\}_{s\in \mathbb{R}}+(2\pi\mathbb{Z}/r_\perp\lambda)^d)\cap
(B_\rho(\alpha_{\xi'})+\{s'\xi'\}_{s'\in \mathbb{R}}+(2\pi\mathbb{Z}/r_\perp\lambda)^d)=\varnothing,$$
for all $\xi,\xi'\in\Lambda,\xi\neq\xi'$.
\el

For each $\xi\in\Lambda$ let us define $A^i_\xi\in \mathbb{S}^{d-1}\cap \mathbb{Q}^d,\ i=1,2,...,d-1$
to be orthogonal vectors to each other and $\xi$. Then for each $\xi\in\Lambda$ we
have that $\{\xi, A^i_\xi,i=1,...,d-1\} \subseteq \mathbb{S}^{d-1}\cap \mathbb{Q}^d$
form an orthonormal basis in $\mathbb{R}^d$. We label by $n_*$ the smallest
natural such that
$$\{n_*\xi, n_*A^i_\xi,i=1,...,d-1\}\subset\mathbb{Z}^d,$$
for every $\xi\in\Lambda$.

Let $\Phi : \mathbb{R}^{d-1} \to \mathbb{R}$ be a smooth function with support in a ball of radius 1. We normalize $\Phi$ such that $\varphi = -\Delta\Phi$ obeys
$$\frac1{(2\pi)^{d-1}}\int_{\mathbb{R}^{d-1}}\varphi^2(x_1,x_2,...,x_{d-1})\dif x_1\dif x_2..\dif x_{d-1}=1.$$
By definition we know $\int_{\mathbb{R}^{d-1}}\varphi \dif x=0$. Define $\psi : \mathbb{R}\to\mathbb{R}$ to be a smooth, mean zero function with
support in the ball of radius 1 satisfying
$$\frac1{2\pi}\int_\mathbb{R}\psi^2(x_d)\dif x_d=1.$$

We define the rescaled cut-off functions
$$\phi_{r_{\perp}}(x_1,x_2,...,x_{d-1})=\frac1{r_{\perp}^{(d-1)/2}}\phi(\frac{x_1}{r_{\perp}},\frac{x_2}{r_{\perp}},...,\frac{x_{d-1}}{r_{\perp}}),$$
$$\Phi_{r_{\perp}}(x_1,x_2,...,x_{d-1})=\frac1{r_{\perp}^{(d-1)/2}}\Phi(\frac{x_1}{r_{\perp}},\frac{x_2}{r_{\perp}},...,\frac{x_{d-1}}{r_{\perp}}),$$
$$\psi_{r_{\parallel}}=\frac1{r_{\parallel}^{1/2}}\psi(\frac{x_d}{r_{\parallel}}).$$
We periodize $\phi_{r_{\perp}},\Phi_{r_{\perp}},\psi_{r_{\parallel}}$ so that they are viewed as periodic functions on $\mathbb{T}^{d-1}, \mathbb{T}^{d-1}$ and $\mathbb{T}$ respectively.
Consider a large real number $\lambda$ such that $\lambda r_{\perp}\in \mathbb{N}$, and a large time oscillation parameter $\mu > 0$. For every $\xi\in\Lambda$ we introduce
$$\psi_{(\xi)}(t,x):=\psi_{\xi,r_{\perp},r_{\parallel},\lambda,\mu}(t,x):=
\psi_{r_{\parallel}}(n_*r_{\perp}\lambda(x\cdot \xi+\mu t)),$$
$$\Phi_{(\xi)}(x):=\Phi_{\xi,r_{\perp},r_{\parallel},\lambda,\mu}(x):=
\Phi_{r_{\perp}}(n_*r_{\perp}\lambda(x-\alpha_\xi)\cdot A^1_{\xi},...,n_*r_{\perp}\lambda(x-\alpha_\xi)\cdot A^{d-1}_{\xi}),$$
$$\phi_{(\xi)}(x):=\phi_{\xi,r_{\perp},r_{\parallel},\lambda,\mu}(x):=
\phi_{r_{\perp}}(n_*r_{\perp}\lambda(x-\alpha_\xi)\cdot A^1_{\xi},...,n_*r_{\perp}\lambda(x-\alpha_\xi)\cdot A^{d-1}_{\xi}),$$
where $\alpha_\xi\in\mathbb{R}^d$ are shifts and $\lambda$ large enough to ensure that $\{\Phi_{(\xi)}\}_{\xi\in\Lambda}$ have mutually disjoint support by Lemma \ref{lem:4}.

The intermittent jets $W_{(\xi)} :\mathbb{R} \times\mathbb{T}^d \to \mathbb{R}^d$ are defined as
$$W_{(\xi)}(t,x):=W_{\xi,r_{\perp},r_{\parallel},\lambda,\mu}(t,x):=\xi\psi_{(\xi)}(t,x)\phi_{(\xi)}(x).$$

By the choice of $\alpha_\xi$ and definition of $W_{(\xi)}$ we have that
\begin{align}
W_{(\xi)}\otimes W_{(\xi')}\equiv 0,\ for\ \xi\neq\xi'\in\Lambda,\label{*}\\
\frac{1}{(2\pi)^d}\int_{\mathbb{T}^d}W_{(\xi)}(t,x)\otimes W_{(\xi)}(t,x)\dif x=\xi\otimes\xi,\label{***}\\
\div(W_{(\xi)}\otimes W_{(\xi)}))=\frac{1}{\mu}\partial_t(\psi^2_{(\xi)}\phi^2_{(\xi)}\xi).\label{*5}
\end{align}

These facts combined with Lemma \ref{lem:3} imply that
$$\frac{1}{(2\pi)^d}\sum_{\xi\in\Lambda}\gamma_{\xi}^2(R)\int_{\mathbb{T}^d}W_{(\xi)}(t,x)\otimes W_{(\xi)}(t,x)\dif x=R,$$
for every symmetric matrix $R$ satisfying $|R-\Id|\leq 1/2$.

Since $W_{(\xi)}$ are not divergence free, inspired by \cite[Section 4.1]{AX} we
introduce the corrector term
\begin{align}
V_{(\xi)}:=\frac{1}{(n_*\lambda)^2}(\xi\otimes\nabla\Phi_{(\xi)}-\nabla\Phi_{(\xi)}\otimes\xi)\psi_{(\xi)},\notag
\end{align}
which is a skew-symmetric matrix.

Then by a direct computation
\begin{align}
\div V_{(\xi)}&=\frac{1}{(n_*\lambda)^2}\psi_{(\xi)}(\xi \Delta\Phi_{(\xi)}-(\xi\cdot\nabla)(\nabla\Phi_{(\xi)}))+\frac{1}{(n_*\lambda)^2}\xi\nabla\Phi_{(\xi)}\cdot\nabla\psi_{(\xi)}-\frac{1}{(n_*\lambda)^2}\nabla\Phi_{(\xi)}\xi\cdot\nabla\psi_{(\xi)}\notag\\
&=\psi_{(\xi)}\phi_{(\xi)}\xi -\frac{1}{(n_*\lambda)^2}\nabla\Phi_{(\xi)}\xi\cdot\nabla\psi_{(\xi)}=W_{(\xi)} -\frac{1}{(n_*\lambda)^2}\nabla\Phi_{(\xi)}\xi\cdot\nabla\psi_{(\xi)}.\label{divOmega}
\end{align}

Finally, we claim that for $N, M \geq 0$ and $p\in [1, \infty]$ the
following holds
\begin{align}
\|\nabla^N\partial_t^M\psi_{(\xi)}\|_{C_tL^p}&\lesssim r_\parallel^{\frac{1}{p}-\frac{1}{2}}(\frac{r_\perp\lambda}{r_\parallel})^N(\frac{r_\perp\lambda\mu}{r_\parallel})^M,\label{int2}\\
\|\nabla^N\phi_{(\xi)}\|_{L^p}+\|\nabla^N\Phi_{(\xi)}\|_{L^p}&\lesssim r_\perp^{\frac{d-1}{p}-\frac{d-1}2}\lambda^N,\label{int3}\\
\|\nabla^N\partial_t^MW_{(\xi)}\|_{C_tL^p}+\lambda\|\nabla^N\partial_t^MV_{(\xi)}\|_{C_tL^p}&\lesssim r_\perp^{\frac{d-1}{p}-\frac{d-1}2}r_\parallel^{\frac{1}{p}-\frac{1}{2}}\lambda^N(\frac{r_\perp\lambda\mu}{r_\parallel})^M,\label{int4}
\end{align}
where implicit constants may depend on $p,N$ and $M$, but are independent of $\lambda,r_\perp,r_\parallel,\mu$. These estimates can be easily deduced from the definitions.


\section{Construction of global probabilistically strong solutions}\label{cogpss}
This section is devoted to the proof of Theorem \ref{thm:4}. The goal is to establish existence of non-unique global-in-time probabilistically strong solutions to the power-law system
(\ref{eq:pl}) for every given divergence free initial condition in $L_\sigma^2\cap W^{1,s}$ for $s\in[r^*,\frac{2d}{d+2})$.

Following the way in \cite{HZZ} we first construct the local probabilistically strong solutions with Cauchy problem. Then
we use the final value at the stopping time of any
such convex integration solution as a new initial condition for the convex integration procedure. This way we are able to extend the convex integration solutions as probabilistically strong solutions defined on the whole time interval $[0, \infty)$.

To this end, we intend to prescribe an arbitrary random initial condition $u_0\in L^2_\sigma\cap W^{1,s}$ $\mathbf{P}$-a.s. independent of the given Wiener process $B$. Let $(\mathcal{F}_t)_{t> 0}$ be the augmented joint canonical filtration on $(\Omega, \mathcal{F})$ generated by $B$ and $u_0$. Then $B$ is a $(\mathcal{F}_t)_{t> 0}$-Wiener process and $u_0$ is $\mathcal{F}_0$-measurable. We include the initial value into the linear part $z$, namely, we let $z$ satisfy the following stochastic equation.
\begin{align}\label{linear}
\dif z-\Delta z\dif t+\nabla\pi_1\dif t&=\dif B,\notag\\
\div z&=0,\\
z(0)&=u_0,\notag
\end{align}
and
\begin{align}\label{eq:v1}
\partial_tv+\div((v+z)\otimes(v+z))-\div\mathcal{A}(Dv+Dz)+\Delta z+\nabla\pi_2&=0,\notag\\
\div v&=0,\notag\\
v(0)&=0.
\end{align}
We consider an increasing sequence $\{\lambda_q\}_{q\in\mathbb{N}_0}\subset\mathbb{N}$ which diverges to $\infty$, and a sequence $\{\delta_q\}_{q\in\mathbb{N}}\subset (0,1]$
 which is decreasing to 0. We choose $a\in\mathbb{N},b\in\mathbb{N},\beta\in(0,1)$ and let
$$\lambda_q=a^{(b^q)},\ \ \delta_q=\lambda_1^{2\beta}\lambda_q^{-2\beta},$$
where $\beta$ will be chosen sufficiently small and $a,b$ will be chosen sufficiently large. In addition, we used $\sum_{q\geq1} \delta_q^{1/2}\leq  \sum_{q\geq1}a^{b\beta-qb\beta}\leq \frac{1}{1-a^{-b\beta}}\leq2$
which boils down to
\begin{align}\label{ab2}
a^{b\beta}\geq2
\end{align}
assumed from now on. More details on the choice of these parameters will be given below in the course of the construction.

The iteration is indexed by a parameter $q\in\mathbb{N}_0$. At each step $q$, a pair $(v_q, \mathring{R}_q)$ is constructed solving the following system
\begin{align}\label{3.1}
\partial_tv_q+\div((v_q+z_q)\otimes(v_q+z_q))-\div\mathcal{A}(Dv_q+Dz_q)+\Delta z_q+\nabla\pi_{q}
&=\div\mathring{R}_q,\notag\\
\div v_q&=0,\\
v_q(0)&=0.\notag
\end{align}
Here we decompose $z=z^{in}+Z$ with $z^{in}(t)=e^{t\Delta}u_0$, where $e^{t\Delta}$ is heat semigroup, and define $z_q=z^{in}+Z_q=z^{in}+\mathbb{P}_{\leq f(q)}Z$ with $f(q)=\lambda_{q+1}^{\frac{\alpha}{4+2d}}$,  $\alpha\in(0,1)$ is given in Section \ref{Cp}. So $u(0)=v(0)+z(0)=\lim_{q\to\infty}v_q(0)+u_0=u_0$.

Let $L,M\geq1$ be given and assume in addition that $\mathbf{P}$-a.s.
\begin{align}
\|u_0\|_{L^2}+\|u_0\|_{W^{1,s}}\leq M.\label{omegan}
\end{align}

Suppose that $\mathrm{Tr}((-\Delta)^{2\sigma}GG^*)<\infty$ for some $\sigma>0$, then for fixed $\delta\in(0,\frac{1}{24})$, the following result
follows from \cite[Theorem 5.14]{DPZ92} together with together with Kolmogorov continuity criterion.
\bp
Suppose that $\mathrm{Tr}((-\Delta)^{2\sigma}GG^*)<\infty$ for some $\sigma>0$. Then for $\delta\in(0,\frac{1}{24})$ and $T>0$
$$\mathbb{E}^{\mathbf{P}}[ \|Z\|_{C_TH^{1+\sigma}}+\|Z\|_{C_T^{1/2-2\delta}L^2} ]<\infty.
$$
\ep

\br\label{jieshiz}
Compared to Navier-Stokes equations, power-law system have $\mathcal{A}(Du+Dz)$ with $\mathcal{A}$ being a nonlinear operator. To give a meaning to this term, it requires $z\in H^1$, which leads to our condition on the noise.
\er

We define
$$T_L:=\inf\{t\geq0,\|Z(t)\|_{H^{1+\sigma}}\geq L/C_S\}\wedge\inf\{t\geq0,\|Z\|_{C_t^{1/2-2\delta}L^2}\geq L/C_S\}\wedge L,$$
where $C_S$ is from the Sobolev embedding $\|f\|_{L^\infty}\leq C_S \|f\|_{H^{\frac{d+\sigma}{2}}}$, then $T_L$ is P-a.s. strictly positive stopping time. Moreover, by calculating, the following terms are controlled by $L\lambda_{q+1}^{\frac\alpha4}$, i.e. for $t\in[0,T_L]
$
\begin{align}\label{3.2}
\|Z_q(t)\|_{L^\infty}, \|\nabla Z_q(t)\|_{L^\infty}, \|\nabla^2 Z_q(t)\|_{L^\infty}, \|Z_q\|_{C_t^{1/2-2\delta} L^{\infty}}, \|\nabla Z_q\|_{C_t^{1/2-2\delta}L^2}\leq L\lambda_{q+1}^{\alpha/4}.
\end{align}

We suppose that there is a deterministic constant
$$M_L\geq\max\{4(L+ M)^2,(2C_{\nu_0}\vee8)(2\pi)^{\frac d2}(L+ M),(2C_{\nu_1}\vee4C_{r,\nu_0,\nu_1})(2\pi)^{\frac{d(r-1)}{2}}(L+ M)^{r-1}\},$$
where constants $C_{\nu_0},C_{\nu_1}$ and $C_{r,\nu_0,\nu_1}$ come from Lemma \ref{lem:1}.

We denote $A=2M_L,\sigma_q=2^{-q},q\in\mathbb{N}_0\cup\{-1\},\gamma_q=2^{-q}
,q\in\mathbb{N}_0\backslash\{3\},\gamma_3=K$. Here $K>1$ is
a large constant which used to distinguish different solutions.

Under the above assumptions, our main iteration reads as follows:

\bp\label{prop:1}
Let $L,M\geq 1$ and assume (\ref{omegan}), there exists a choice of parameters $a,b,\beta$ such
that the following holds true: Let $(v_q, \mathring{R}_q)$ be an $(\mathcal{F}_t)_{t\geq0}$-adapted solution to (\ref{3.1}) satisfying
\begin{align}\label{3.3}
\|v_q(t)\|_{L^2}\leq\begin{cases}
M_0(M_L^{1/2}\sum_{m=1}^q\delta_m^{1/2}+\sum_{m=1}^q\gamma_m^{1/2})+3M_0(M_L+qA)^{1/2},&t\in(\frac{\sigma_{q-1}}{2}\wedge T_L,T_L],\\
0,&t\in[0,\frac{\sigma_{q-1}}{2}\wedge T_L],
\end{cases}
\end{align}
for a universal constant $M_0$, and
\begin{align}
\|v_q\|_{C_{t,x}^1}&\leq\lambda_q^{d+1}M_L^{1/2},\ \ t\in[0,T_L],\label{3.4}\\
\|v_q\|_{C_{t,x}^2}&\leq\lambda_q^{\frac{3d}{2}+2}M_L^{1/2},\ \ t\in[0,T_L],\label{3.5}\\
\|\mathring{R}_q(t)\|_{L^1}&\leq\delta_{q+1}M_L,\ \ t\in(\sigma_{q-1}\wedge T_L,T_L],\label{3.6}\\
\|\mathring{R}_q(t)\|_{L^1}&\leq M_L+qA,\ \ t\in[0,T_L].\label{3.7}
\end{align}
Then there exists an $(\mathcal{F}_t)_{t>0}$-adapted process $(v_{q+1}, \mathring{R}_{q+1})$ which solves (\ref{3.1}) and satisfies
\begin{align}
\|v_{q+1}-v_q\|_{L^2}\leq
\begin{cases}
M_0(M_L^{1/2}\delta_{q+1}^{1/2}+\gamma_{q+1}^{1/2}),&t \in(4\sigma_q\wedge T_L,T_L],\\
M_0((M_L+qA)^{1/2}+\gamma_{q+1}^{1/2}),\ &t\in(\frac{\sigma_q}2\wedge T_L,4\sigma_q\wedge T_L],\\
0,&t\in[0,\frac{\sigma_q}2\wedge T_L],
\end{cases}\label{3.8}\\
\|\mathring{R}_{q+1}(t)\|_{L^1}\leq
\begin{cases}
M_L\delta_{q+2},&t \in(\sigma_q\wedge T_L,T_L],\\
M_L\delta_{q+2}+\sup_{s\in[t-\sigma_q,t]}\|\mathring{R}_{q}(s)\|_{L^1},&t\in(\frac{\sigma_q}2\wedge T_L,\sigma_q\wedge T_L],\\
\sup_{s\in[t-\sigma_q,t]}\|\mathring{R}_{q}(s)\|_{L^1}+A,&t\in[0,\frac{\sigma_q}2\wedge T_L],
\end{cases}\label{3.9}
\end{align}
and $(v_{q+1}, \mathring{R}_{q+1})$ obeys (\ref{3.3})-(\ref{3.7}) at the level $q+1$. Furthermore,
\begin{align}
\|v_{q+1}(t)-v_q(t)\|_{W^{1,s}}\leq M_L^{1/2}\delta_{q+1}^{1/2},\ \ t\in[0,T_L],\label{3.10}
\end{align}
and for $t\in(4\sigma_q\wedge T_L,T_L]$ we have
\begin{align}
|\|v_{q+1}\|_{L^2}^2-\|v_{q}\|_{L^2}^2-d\gamma_{q+1}|\leq(2d+1)M_L\delta_{q+1}.\label{3.11}
\end{align}
\ep

The proof of this result is presented in Section \ref{proof:prop1} below.

\br\label{BMSno}
Comparing to \cite[Theorem C]{BMS}, we only obtain solutions in $L^p(0, T_L;L^2)$ for every $p \in[1, \infty)$ but not for $p = \infty$. We refer to \cite[Remark 5.3]{HZZ} for explanation on this point.
\er

We intend to start the iteration from $v_0\equiv0$ on $[0, T_L]$. In that case by (\ref{rdeltav}) we have $\mathring{R}_0=z_0\mathring{\otimes}z_0-\mathcal{A}(Dz_0)+2D{z_0}$, which by Lemma \ref{lem:1} and H\"older's inequality implies that
\begin{align}
\|\mathcal{A}(Dz_0)\|_{L^1}\leq
\begin{cases}
C_{\nu_1}\|Dz_0\|_{L^1}^{r-1},&for\ \nu_0=0,r\leq2,\\
C_{\nu_0}\|Dz_0\|_{L^1},&for\ \nu_0>0,r\leq2,\\
C_{r,\nu_0,\nu_1}\|Dz_0\|_{L^{r-1}}(1+\|Dz_0\|^{r-2}_{L^{r-1}}),&for\ r> 2.
\end{cases}\notag
\end{align}
Moreover, by
$\|Dz_0\|_{L^{r^*}}\leq (\|u_0\|_{W^{1,s}}+\|Z\|_{H^1})(2\pi)^{\frac{d}{2}}\leq (2\pi)^{\frac{d}{2}}(L+M
)$, and the choice of $M_L$ we have
\begin{align}\label{adz0}
\|\mathcal{A}(Dz_0)\|_{L^1}\leq \frac{1}{2}M_L,
\end{align}
which yields that
\begin{align}
\|\mathring{R}_0(t)\|_{L^1}&\leq \|z_0\|_{L^{2}}^2+\frac{1}{2}M_L+2\|\nabla z_0\|_{L^{1}}\notag\\
&\leq (\|z^{in}\|_{L^2}+\|Z_0\|_{L^2})^2+\frac{1}{2}M_L+2\|z^{in}\|_{W^{1,1}}+2\|Z_0\|_{W^{1,1}}\notag\\
&\leq (M+L)^2+\frac{1}{2}M_L+2(2\pi)^{\frac{d}{2}}(M+L) \leq M_L.
\end{align}
So (\ref{3.6}) as well as (\ref{3.7}) are satisfied on the level $q = 0$, since $\delta_1 = 1.$

We deduce the following result.
\bt\label{thm:3}
Let $r\in(1,\frac{3d+2}{d+2}),s\in[r^*,\frac{2d}{d+2})$. Then there exists a $\mathbf{P}$-a.s. strictly positive stopping time $T_L$, arbitrarily large by choosing $L$ large, such that for any initial condition $u_0\in L^2_\sigma\cap W^{1,s}\ \mathbf{P}$-a.s. independent of the Brownian motion $B$ the following holds true: There exists an $(\mathcal{F}_t)_{t> 0}$-adapted process $u$ which belongs to $L^p([0,T_L];L^2)\cap C([0, T_L];W^{1,
s})\ \mathbf{P}$-a.s. for all $p \in[1, \infty)$, and is an analytically weak solution to (\ref{eq:pl}) with $u(0) = u_0$. There are infinitely many such solutions u.
\et

\begin{proof}
Let the additional assumption (\ref{omegan}) be satisfied for some $M > 1$. Letting $v_0\equiv0$ we
repeatedly apply Proposition \ref{prop:1} and obtain  $(\mathcal{F}_t)_{t\geq0}$-adapted processes $(v_q,\mathring{R}_q)$, we obtain $v_q\to v$ in $C([0, T_L];W^{1,s})$ as a consequence of (\ref{3.10}). Similar argument as \cite[Theorem 5.4]{HZZ} and using (\ref{3.8}) and (\ref{3.9}) we obtain $v_q\to v$ in $L^p([0,T_L];L^2)$ and $\mathring{R}_q\to0$ in $L^p([0,T_L];L^1)$ for all $p\in[1,\infty)$.
Moreover, we have
$$|\int_0^t\int[\mathcal{A}(Dv_q+Dz_q)-\mathcal{A}(Dv+Dz)]\nabla\varphi|\dif x\dif s\leq
\|\nabla\varphi\|_{L_t^1L^{\infty}}\|\mathcal{A}(Dv_q+Dz_q)-\mathcal{A}(Dv+Dz)\|_{C_tL^1}.$$
By Lemma \ref{lem:1} and H\"older's inequality we obtain for $t\in[0,T_L]$ and
for $r<2,\nu_0>0$
\begin{align*}
\|\mathcal{A}(Dv_q+Dz_q)-\mathcal{A}(Dv+Dz)\|_{C_tL^1}&\lesssim(\|Dv_q-Dv\|_{C_tL^{1}}+\|Dz_q-Dz\|_{C_tL^{1}})\\
&\lesssim\|v_q-v\|_{C_tW^{1,1}}+\|{z}_{q}-z\|_{C_tH^1}\to0,
\end{align*}
and for $r<2,\nu_0=0$
\begin{align*}
\|\mathcal{A}(Dv_q+Dz_q)-\mathcal{A}(Dv+Dz)\|_{C_tL^1}\lesssim\|v_{q}-v\|_{C_tW^{1,1}}^{r-1}+\|{z}_{q}-z\|_{C_tH^1}^{r-1}\to0,
\end{align*}
and for $r\geq2$,
\begin{align*}
\|\mathcal{A}(Dv_q+Dz_q)&-\mathcal{A}(Dv+Dz)\|_{C_tL^1}\\
&\lesssim(\|Dv_q-Dv\|_{C_tL^{r-1}}+\|Dz_n-Dz\|_{C_tL^{r-1}})(1+\|Dv\|_{C_tL^{r-1}}^{r-2}+\|Dz\|_{C_tL^{r-1}}^{r-2}).\\
&\lesssim(\|v_{q}-v\|_{C_tW^{1,r-1}}+\|{z}_{q}-z\|_{C_tH^1})(1+\|v\|_{C_tW^{1,r-1}}^{r-2}+\| z\|_{C_tH^1}^{r-2})\to 0.
\end{align*}

Finally, we let $q\to\infty$ in (\ref{3.1}) and obtain the process $u = v + z$ satisfies (\ref{eq:pl}) before $T_L$ in the analytic weak sense. Since $v_q(0) = 0$ we deduce $v(0) = 0$, which implies that $u(0) = u_0$.

By a similar argument as \cite[Theorem 5.4]{HZZ} we obtain nonuniqueness of the constructed solutions on $[0,T_L]$ and the result holds for general initial condition $u_0\in L_\sigma^2\cap W^{1,s}$.

For a general initial condition $u_0\in L_\sigma^2\cap W^{1,s},\ \mathbf{P}$-a.s., define $\Omega_M:=\{ M-1\leq\|u_0\|_{L^2}+\|u_0\|_{W^{1,s}}< M\}$. Then the first part of this proof gives the existence of infinitely many adapted solutions $u^M$ on each $\Omega_M$.
Letting $u:=\sum_{M\in \mathbb{N}}u^M1_{\Omega_M}$ concludes the proof.
\end{proof}

Now we have all in hand to complete the proof of  Theorem \ref{thm:4}.
\begin{proof}[Proof of Theorem \ref{thm:4}]
By Theorem \ref{thm:3} we constructed a probabilistically strong solution $u$ before
the stopping time $T_L$ starting from the given initial condition $u_0\in L^2_\sigma\cap W^{1,s}\ \mathbf{P}$-a.s. By a similar argument as \cite[Theorem 1.1]{HZZ} we obtain $\|u(T_L)\|_{L^2}<\infty\ \mathbf{P}$-a.s.. Moreover,
\begin{align}
\|u(T_L)\|_{W^{1,s}}&\leq\|z(T_L)\|_{W^{1,s}}+\sum_{q\geq0}\|v_{q+1}(T_L)-v_q(T_L)\|_{W^{1,s}}\notag\lesssim \sum_{q\geq0}M_L^{1/2}\delta_{q}^{1/2}<\infty.\notag
\end{align}
This implies that $\|u(T_L)\|_{L^2}+\|u(T_L)\|_{W^{1,s}} < \infty\ \mathbf{P}$-a.s. hence we can use the value $u(T_L)$ as a new initial condition in Theorem \ref{thm:3}. Then a similar argument as \cite[Theorem 1.1]{HZZ} we obtain the result.
\end{proof}

\begin{proof}[Proof of Corollary \ref{thm:10}]
The proof follows from the same argument as \cite[Corollary 1.2]{HZZ}.
\end{proof}

\subsection{Proof of Proposition \ref{prop:1}}\label{proof:prop1}
The proof proceeds in several main steps which are the same in many convex integration schemes. First of all, we start the construction by fixing the parameters in Section~\ref{Cp} and proceed with a mollification step in Section~\ref{subsub2}. Section~\ref{subsub3} introduces the new iteration $v_{q+1}$. This is the main part of the construction which used the generalized intermittent jets introduced in Section \ref{gij}.  Section~\ref{subsub4} contains the inductive estimates of $v_{q+1}$, whereas in Section~\ref{subsub5} we show how the energy is controlled. Finally, in Section~\ref{cotrs}, we define the new stress $\mathring{R}_{q+1}$ and  establish the inductive moment estimate on $\mathring{R}_{q+1}$ in Section~\ref{subsub7}.
\subsubsection{Choice of parameters}\label{Cp}
In the sequel, additional parameters will be indispensable and their value has to be carefully chosen in order to respect all the compatibility conditions appearing in the estimates below. First, for a sufficiently small $\alpha\in(0,1)$ to be chosen, we take $l:=\lambda_{q+1}^{-\frac{3\alpha}{2}}\lambda_q^{-2}$ and have $l^{-1}\leq\lambda_{q+1}^{2\alpha}$ provided $\alpha b>4$. For the choice of $l$ we require
\begin{align}
l^{1/6}\leq\sigma_q^{d/4+1/2},\label{a1}
\end{align}
which is satisfied if $a$ large enough.

Now we introduce $N=[\frac{(2d+1)(2-s)}{(2-s)d-2s}]+d+3$. In the sequel, we also need
$$(2d+8)(\frac{[dN]}2+2)\alpha<\frac{1}{N},\ \ \alpha b>\max\{16(d+1),\frac{2(d+1)(2d+4)}{\sigma}\},$$ $$\alpha>\max\{32, \frac{4(2d+4)}{\sigma}\}\beta b,\ \
\frac1{2N}-(12d+27)\alpha>2\beta b,\ \  P(s)-(10d+19)\alpha>2\beta b.$$
 In the case $\nu_0=0$ we additionally require that
$$\alpha> \frac{16\beta b}{r-1},\ \ \alpha>\frac{2d+4}{\sigma(r-1)}2\beta b,\ \ P(s)-(10d+17)\alpha> \frac{16\beta b}{r-1}.$$
Here we denote $P(x)=(\frac{1}{x}-\frac{1}{2})(d-\frac{2d+1}{N})-1$, and by the definition of $N$ we obtain $P(s)>0$.

The above can be obtained by choosing $\alpha$ small such that $\frac1{2N}-(12d+27)\alpha>\alpha$,  $P(s)-(10d+19)\alpha>\alpha$ and $(2d+8)(\frac{[dN]}2+2)\alpha<\frac{1}{N}$, and choosing $b\in N\mathbb{N}$ large enough such that $b>\frac{1}{\alpha}\max\{16d+16,\frac{2(d+1)(2d+4)}{\sigma}\}$, and finally choosing $\beta$ small such that $\alpha>\beta b\max\{\frac{16}{r-1},32,\frac{2(2d+4)}{\sigma(r-1)},\frac{4(2d+4)}{\sigma}\}$.

From the above we obtained
\begin{align}
M_L + K + qA\leq l^{-1},\ \ (M_L+K+qA)(\lambda_{q+1}^{-\alpha/8+2\beta b}+\lambda_{q+1}^{-\frac{\alpha\sigma}{2d+4}+2\beta b })\ll1, \label{a2}
\end{align}
by choosing $a$ large enough.

In the sequel, we increase $a$ in order to absorb various implicit and universal constants in the following estimates.
So the last free parameter $a$ is given through (\ref{ab2}), (\ref{a1}) and (\ref{a2}).

\subsubsection{Mollification}\label{subsub2}
We replace $v_q$ by a mollified field $v_l$, and we define
\begin{align}
v_l=(v_q*_x\phi_l)*_t\varphi_l,\ \ \mathring{R}_l=(\mathring{R}_q*_x\phi_l)*_t\varphi_l\notag,\ \ z_l=(z_q*_x\phi_l)*_t\varphi_l,
\end{align}
where $\phi_l:=\frac{1}{l^d}\phi(\frac{\cdot}{l})$ is a family of standard mollifiers on $\mathbb{R}^d$, and $\varphi_l:=\frac{1}{l}\varphi(\frac{\cdot}{l})$ is a family of standard mollifiers with support in $(0,1)$. The one side mollifier is used to preserve adaptedness, and we extend $v_q,\mathring{R}_q$ and $z_q$ to $t \leq0$ by taking them equal to the value at $t = 0$. Then $(v_q,z_q,\mathring{R}_q)$ also satisfies equation (\ref{3.1}) for $t<0$ as $\partial_tv_q(0)=0$ from our construction. It is easy to see that $v_l,\mathring{R}_l$ and $z_l$ are $(\mathcal{F}_t)_{t\geq0}$-adapted.

By calculating and (\ref{3.1}), it follows that $(v_l,\mathring{R}_l)$  satisfies
\begin{align}
\partial_tv_l+\div((v_l+z_l)\otimes(v_l+z_l))-\div\mathcal{A}(Dv_l+Dz_l)+\Delta z_l+\nabla\pi_{l}
=\div(\mathring{R}_l+R_{com1}+{R}_{com2}),\label{eq:v_l}
\end{align}
where
\begin{align}
R_{com1}&=(v_l+z_l)\mathring{\otimes}(v_l+z_l)-((v_q+z_q)\mathring{\otimes}(v_q+z_q))*_x\phi_l*_t\varphi_l,\notag \\
R_{com2}&=(\mathcal{A}(Dv_q+Dz_q)*_x\phi_l)*_t\varphi_l-\mathcal{A}(Dv_l+Dz_l),\notag\\
\pi_l&=(\pi_q*_x\phi_l)*_t\varphi_l-\frac{1}{d}\tr R_{com1}.\notag
\end{align}
Here $\mathring{R}$ means the trace-free part of $R$. We used the nonlinear term $\mathcal{A}(Q)$ is  trace-free when $Q$ is trace-free and $Dv,Dz$ is trace-free due to divergence free condition.

\subsubsection{Construction of $v_{q+1}$}\label{subsub3}
Let us now proceed with the construction of the perturbation $w_{q+1}$ which then defines the next iteration by $v_{q+1} := v_l +w_{q+1}$. To this end, we employ the generalized intermittent jets introduced in section \ref{gij}. We choose the following parameters
\begin{align}\label{eq:sym}
\lambda=\lambda_{q+1},\ r_{\parallel}=\lambda_{q+1}^{\frac{d+2}N-1},\ r_{\perp}=\lambda_{q+1}^{\frac1N-1},\ \mu=\lambda_{q+1}^{\frac{d}2},
\end{align}
where we recall $N=[\frac{(2d+1)(2-s)}{(2-s)d-2s}]+d+3$. It is required that $b$ is a multiple of $N$ to ensure that $\lambda_{q+1}r_{\perp}=a^{(b^{q+1})/N}\in\mathbb{N}$.\\
\begin{remark}\label{rem:1}
To explain how to choose $r_{\parallel}, r_{\perp}, \mu$, we set for $\alpha_0,\beta_0\in(0,1),c>0$
$$ r_{\parallel}=\lambda_{q+1}^{\alpha_0-1},\ r_{\perp}=\lambda_{q+1}^{\beta_0-1},\ \mu=\lambda_{q+1}^{c}.$$
First note that $\lambda_{q+1}^{-1}\ll r_\perp\ll r_\parallel\ll 1$ implies $0<\beta_0<\alpha_0<1$.
The estimates (\ref{2.18*}), (\ref{2.27*}) and (\ref{rlin1}) below require  $ r_\parallel^{-\frac{1}{2}}r_\perp^{-\frac{d-1}{2}}\mu^{-1},\lambda r_\perp^{\frac{d-1}{s}-\frac{d-1}{2}} r_\parallel^{\frac{1}{s}-\frac{1}{2}},r_\perp^{\frac{d+1}{2}} r_\parallel^{-\frac{1}{2}}\mu$ small,  i.e.
$$\frac d2-\frac{1}{2}\alpha_0-\frac {d-1}2\beta_0<c<\frac{d}{2}+\frac{1}{2}\alpha_0-\frac {d+1}2\beta_0,$$ $$\frac{1}{s}-\frac{1}{2}>\frac{1}{d-(\alpha_0+(d-1)\beta_0)}.$$
By $0<\beta_0<\alpha_0$, the last inequality holds only when $s<\frac{2d}{2+d}$. When $s<\frac{2d}{2+d}$ we set $\beta_0=\frac{1}{N},\alpha_0=\frac{d+2}{N},c=\frac{d}2$ and the last inequality holds if $N>\frac{(2d+1)(2-s)}{(2-s)d-2s}$. We also require $\alpha_0<1$ i.e. $N>d+2$, so we finally choose $N=[\frac{(2d+1)(2-s)}{(2-s)d-2s}]+d+3$.
\end{remark}
As the next step, we shall define certain amplitude functions used in the definition of the perturbations $w_{q+1}$. Similarly as \cite[Section 5.1.2]{HZZ} we define
\begin{align}
\rho:=2\sqrt{l^2+|\mathring{R}_l|^2}+\frac{\gamma_{q+1}}{(2\pi)^d}.
\end{align}
It follows that $\rho$ is $(\mathcal{F}_ t)_{t> 0}$-adapted. Then by the estimate of $\rho$ in Appendix \ref{s:appA} we have for $N\geq1$
$$\|\rho\|_{C_{[4\sigma_q\wedge T_L,T_L],x}^N}\lesssim l^{2-(d+4)N}\delta_{q+1}M_L+\gamma_{q+1},$$
and for $t \in [0, T_L]$
$$\|\rho\|_{C_{t,x}^N}\lesssim l^{2-(d+4)N}(M_L+qA)+\gamma_{q+1}.$$

We define the amplitude functions
\begin{align}
a_{(\xi)}(\omega,x,t):=\rho^{1/2}(\omega,x,t)\gamma_{\xi}(\Id-\frac{\mathring{R}_l(\omega,x,t)}{\rho(\omega,x,t)}), \label{3.13*}
\end{align}
where $\gamma_\xi$ is introduced in Lemma \ref{lem:1}. Since
we have
$$|\Id-\frac{\mathring{R}_l}{\rho}-\Id|\leq\frac{1}{2},$$
by Lemma \ref{lem:3}, it holds
\begin{align}
\rho\ \Id-\mathring{R}_l=\sum_{\xi\in\Lambda}\rho\gamma_{\xi}^2(\Id-\frac{\mathring{R}_l}{\rho})\xi\otimes\xi=\sum_{\xi\in\Lambda}a_{(\xi)}^2\xi\otimes\xi.\label{**}
\end{align}

By the estimate of $a_{(\xi)}$ in Appendix \ref{s:appA} we have for $t\in[4\sigma_q\wedge T_L,T_L]$
\begin{align}
\|a_{(\xi)}(t)\|_{L^2}\leq \frac{M^*}{4|\Lambda|}(M_L^{1/2}\delta_{q+1}^{1/2}+\gamma_{q+1}^{1/2}),\label{al21}\\
\|a_{(\xi)}\|_{C_tL^2}\leq \frac{M^*}{4|\Lambda|}((M_L+qA)^{1/2}+\gamma_{q+1}^{1/2})\ \ for\ t\in[0,T_L]\label{al22}.
\end{align}
And for $N\in\mathbb{N}_0$
\begin{align}
\|a_{(\xi)}\|_{C_{[4\sigma_q\wedge T_L,T_L],x}^N}\lesssim l^{-2d-2-(d+4)N}(\delta_{q+1}^{1/2}M_L^{1/2}+\gamma_{q+1}^{1/2}),
\label{3.12}\\
\|a_{(\xi)}\|_{C_{t,x}^N}\lesssim l^{-2d-2-(d+4)N}((M_L+qA)^{1/2}+\gamma_{q+1}^{1/2}),\ \ for\ t\in[0,T_L].\label{3.13}
\end{align}

With these preparations in hand, we define the principal part
\begin{align}
w_{q+1}^{(p)}:=\sum_{\xi\in\Lambda}a_{(\xi)}W_{(\xi)},\notag
\end{align}
and by (\ref{*}), (\ref{***}) and (\ref{**}) we have
\begin{align}
w_{q+1}^{(p)}\otimes w_{q+1}^{(p)}+\mathring{R}_l&=\sum_{\xi\in\Lambda}a_{(\xi)}^2
(W_{(\xi)}\otimes W_{(\xi)}-\xi\otimes\xi)+\rho \Id\notag\\
&=\sum_{\xi\in\Lambda}a_{(\xi)}^2\mathbb{P}_{\neq0}(W_{(\xi)}\otimes W_{(\xi)})+\rho \Id,\label{2.12}
\end{align}
where we use the notation $\mathbb{P}_{\neq0}f:=f-\frac{1}{(2\pi)^d}\int_{\mathbb{T}^d}f(x)\dif x$.

We define the incompressibility corrector by
\begin{align}
w_{q+1}^{(c)}&:=\sum_{\xi\in\Lambda}-\frac{1}{(n_*\lambda_{q+1})^2}a_{(\xi)}\nabla\Phi_{(\xi)}\xi\cdot\nabla\psi_{(\xi)}+\nabla a_{(\xi)}:V_{(\xi)}.\label{2.13}
\end{align}
Here $\nabla a_{(\xi)}:V_{(\xi)}=\sum_j\partial_ja_{(\xi)}V_{(\xi)}^{ij}$.

By (\ref{divOmega}) we have
\begin{align}
w_{q+1}^{(p)}+w_{q+1}^{(c)}=&\sum_{\xi\in\Lambda}a_{(\xi)}(W_{(\xi)}-\frac{1}{(n_*\lambda_{q+1})^2}\nabla\Phi_{(\xi)}\xi\cdot\nabla\psi_{(\xi)})+\nabla a_{(\xi)}:V_{(\xi)}\notag\\
&=\sum_{\xi\in\Lambda}(a_{(\xi)}\div V_{(\xi)}+\nabla a_{(\xi)}:V_{(\xi)})=\sum_{\xi\in\Lambda}\div (a_{(\xi)}V_{(\xi)}).\label{****}
\end{align}
Since $a_{(\xi)}V_{(\xi)}$ is skew-symmetric we obtain
$$\div ( w_{q+1}^{(p)}+w_{q+1}^{(c)})=0.$$

Next we introduce the temporal corrector
\begin{align}
w_{q+1}^{(t)}&:=-\frac1{\mu}\sum_{\xi\in\Lambda}\mathbb{P}_{\neq0}\mathbb{P}_{H}(a_{(\xi)}^2\psi^2_{(\xi)}\phi^2_{(\xi)}\xi),\label{2.14}
\end{align}
where $\mathbb{P}_H$ is the Helmholtz projection, by a direct computation and (\ref{*5}) we obtain
\begin{align}
&\partial_t w_{q+1}^{(t)}+\sum_{\xi\in\Lambda}\mathbb{P}_{\neq0}(a_{(\xi)}^2 \div(W_{(\xi)}\otimes W_{(\xi)}))\notag
\\
&\qquad= -\frac{1}{\mu}\sum_{\xi\in\Lambda}\mathbb{P}_H\mathbb{P}_{\neq0}\partial_t(a_{(\xi)}^2\psi^2_{(\xi)}\phi^2_{(\xi)}\xi)
+\frac{1}{\mu}\sum_{\xi\in\Lambda}\mathbb{P}_{\neq0}( a^2_{(\xi)}\partial_t(\psi^2_{(\xi)}\phi^2_{(\xi)}\xi))\notag
\\&\qquad= (\Id-\mathbb{P}_H)\frac{1}{\mu}\sum_{\xi\in\Lambda}\mathbb{P}_{\neq0}\partial_t(a_{(\xi)}^2\psi^2_{(\xi)}\phi^2_{(\xi)}\xi)
-\frac{1}{\mu}\sum_{\xi\in\Lambda}\mathbb{P}_{\neq0}((\partial_t a^2_{(\xi)})\psi^2_{(\xi)}\phi^2_{(\xi)}\xi).\label{int1}
\end{align}
Note that the first term on the right hand side can be viewed as a pressure term.

Since $\rho$ and $\mathring{R}_l$ are $(\mathcal{F}_t)_{t>0}$-adapted, we know $a_{(\xi)}$ is $(\mathcal{F}_t)_{t>0}$-adapted. Also $W_{(\xi)},V_{(\xi)},\psi_{(\xi)},\phi_{(\xi)}$ are deterministic, so $w_{q+1}^{(p)},w_{q+1}^{(c)},w_{q+1}^{(t)}$ are also $(\mathcal{F}_t)_{t>0}$-adapted.

Let us introduce a smooth cut-off function
\begin{align*}
\chi(t)=\begin{cases}
0,& t\leq \frac{\sigma_q}2,\\
\in (0,1),& t\in (\frac{\sigma_q}2,{\sigma_q} ),\\
1,&t\geq {\sigma_q}.
\end{cases}
\end{align*}
Note that $\|\chi(t)\|_{C_t}\leq 2^{q+1}$ which has to be taken into account in the estimates of $C_{t,x}^i(i=1,2)$ below.

We define the perturbations $\tilde{w}_{q+1}^{(p)},\tilde{w}_{q+1}^{(c)},\tilde{w}_{q+1}^{(t)}$ as follows:
$$\tilde{w}_{q+1}^{(p)}=w_{q+1}^{(p)}\chi,\ \tilde{w}_{q+1}^{(c)}=w_{q+1}^{(c)}\chi,\ \tilde{w}_{q+1}^{(t)}=w_{q+1}^{(t)}\chi^2.$$
Also $\tilde{w}_{q+1}^{(p)},\tilde{w}_{q+1}^{(c)},\tilde{w}_{q+1}^{(t)}$ are $(\mathcal{F}_t)_{t>0}$-adapted.

Finally, the total perturbation $w_{q+1}$ is defined by
$$w_{q+1}:=\tilde{w}_{q+1}^{(p)}+\tilde{w}_{q+1}^{(c)}+\tilde{w}_{q+1}^{(t)},$$
which is $(\mathcal{F}_t)_{t>0}$-adapted, mean zero and divergence-free.

The new velocity $v_{q+1}$ is defined as
\begin{align}
v_{q+1}=v_l+w_{q+1},\notag
\end{align}
which is  $(\mathcal{F}_t)_{t>0}$-adapted, mean zero and divergence-free.

\subsubsection{Estimates of $w_{q+1}$}\label{subsub4}
First we estimate $L^p$-norm for $p\in(1,\infty)$ and verify the inductive estimates (\ref{3.8}) and (\ref{3.3}) at the level $q+1$.

We first consider bound $\tilde{w}^{(p)}_{q+1}$ in $L^2$ by applying Theorem \ref{ihiot}. By (\ref{al21})
we obtain for $t\in(4\sigma _q\wedge T_L,T_L]$ and for some universal constant $M_0\geq1$
\begin{align}
\|\tilde{w}_{q+1}^{(p)}(t)\|_{L^2}&\lesssim\sum_{\xi\in\Lambda}\|a_{(\xi)}(t)W_{(\xi)}(t)\|_{L^2}\lesssim\sum_{\xi\in\Lambda}(\|a_{(\xi)}(t)\|_{L^2}+(r_\perp\lambda)^{-1/2}\|a_{(\xi)}(t)\|_{C^1})\|W_{(\xi)}(t)\|_{L^2}\notag\\
&\lesssim\sum_{\xi\in\Lambda}(\frac{M}{4|\Lambda|}+\lambda_{q+1}^{(6d+12)\alpha-\frac{1}{2N}})(\delta_{q+1}^{1/2}M_L^{1/2}+\gamma_{q+1}^{1/2})\|W_{(\xi)}\|_{C_tL^2}\leq\frac{M_0}2(\delta_{q+1}^{1/2}M_L^{1/2}+\gamma_{q+1}^{1/2}),\label{2.16}
\end{align}
where we used the conditions on the parameters to deduce $(6d+12)\alpha-\frac{1}{2N}<0$.

And similarly for $t\in(\frac{\sigma_q}{2}\wedge T_L,4\sigma _q\wedge T_L]$
\begin{align}
\|\tilde{w}_{q+1}^{(p)}(t)\|_{L^2}\leq \frac{M_0}2((M_L+qA)^{1/2}+\gamma_{q+1}^{1/2}).\label{2.16''}
\end{align}

For general $L^p$-norm, by (\ref{int2})-(\ref{int4}), (\ref{3.12}), for $t\in(4\sigma _q\wedge T_L,T_L]$ we obtain
\begin{align}
\|w_{q+1}^{(p)}(t)\|_{L^p}&\lesssim \sum_{\xi\in\Lambda}\|a_{(\xi)}(t)\|_{L^\infty}\|W_{(\xi)}(t)\|_{L^p}\lesssim (\delta_{q+1}^{1/2}M_L^{1/2}+\gamma_{q+1}^{1/2})l^{-2d-2}r_\perp^{\frac{d-1}{p}-\frac{d-1}{2}} r_\parallel^{\frac{1}{p}-\frac{1}{2}}\notag\\
&\lesssim(\delta_{q+1}^{1/2}M_L^{1/2}+\gamma_{q+1}^{1/2})\lambda_{q+1}^{(4d+4)\alpha+(\frac{1}{p}-\frac{1}{2})(\frac{2d+1}{N}-d)},\label{2.17}
\end{align}

\begin{align}
\|w_{q+1}^{(c)}(t)\|_{L^p}&\lesssim\sum_{\xi\in\Lambda}\frac{1}{\lambda_{q+1}^2}\|a_{(\xi)}(t)\nabla\Phi_{(\xi)}(t)\xi\cdot\nabla\psi_{(\xi)}(t)\|_{L^p}+\sum_{\xi\in\Lambda}\|\nabla a_{(\xi)}(t):V_{(\xi)}(t)\|_{L^p}\notag\\
&\lesssim (\delta_{q+1}^{1/2}M_L^{1/2}+\gamma_{q+1}^{1/2})r_\perp^{\frac{d-1}{p}-\frac{d-1}{2}} r_\parallel^{\frac{1}{p}-\frac{1}{2}}(l^{-2d-2}\frac{r_\perp}{ r_\parallel}+l^{-3d-6}\frac{1}{\lambda_{q+1}})\notag\\
&\lesssim(\delta_{q+1}^{1/2}M_L^{1/2}+\gamma_{q+1}^{1/2})\lambda_{q+1}^{(6d+12)\alpha+(\frac{1}{p}-\frac{1}{2})(\frac{2d+1}{N}-d)-\frac{d+1}{N}},\label{2.18}
\end{align}

and
\begin{align}
\|w_{q+1}^{(t)}(t)\|_{L^p}&\lesssim\mu^{-1}\sum_{\xi\in\Lambda}\|a_{(\xi)}(t)\|_{L^\infty}^2\|\psi_{(\xi)}(t)\|_{L^{2p}}^2\|\phi_{(\xi)}(t)\|^2_{L^{2p}}\notag\\
&\lesssim (\delta_{q+1}^{1/2}M_L^{1/2}+\gamma_{q+1}^{1/2})^2l^{-4d-4}\mu^{-1}r_\perp^{\frac{d-1}{p}-d+1}r_\parallel^{\frac{1}{p}-1}\notag\\
&\lesssim(\delta_{q+1}^{1/2}M_L^{1/2}+\gamma_{q+1}^{1/2})^2\lambda_{q+1}^{(8d+8)\alpha+(\frac{1}{p}-1)(\frac{2d+1}{N}-d)-\frac{d}{2}}.\label{2.18*}
\end{align}

Similarly we obtain for $t\in(\frac{\sigma_q}{2}\wedge T_L,4\sigma _q\wedge T_L]$

\begin{align}
\|w_{q+1}^{(p)}(t)\|_{L^p}
&\lesssim((M_L+qA)^{1/2}+\gamma_{q+1}^{1/2})\lambda_{q+1}^{(4d+4)\alpha+(\frac{1}{p}-\frac{1}{2})(\frac{2d+1}{N}-d)},\label{2.18+-}\\
\|w_{q+1}^{(c)}(t)\|_{L^p}
&\lesssim((M_L+qA)^{1/2}+\gamma_{q+1}^{1/2})\lambda_{q+1}^{(6d+12)\alpha+(\frac{1}{p}-\frac{1}{2})(\frac{2d+1}{N}-d)-\frac{d+1}{N}},\label{2.18--}
\end{align}

and
\begin{align}
\|w_{q+1}^{(t)}(t)\|_{L^p}
&\lesssim((M_L+qA)^{1/2}+\gamma_{q+1}^{1/2})^2\lambda_{q+1}^{(8d+8)\alpha+(\frac{1}{p}-1)(\frac{2d+1}{N}-d)-\frac{d}{2}}.\label{2.18**}
\end{align}

Combining (\ref{2.16}), (\ref{2.18}) and (\ref{2.18*}) we obtain for $t\in(4\sigma _q\wedge T_L,T_L]$
\begin{align}
\|w_{q+1}\|_{C_tL^2}&\leq (\delta_{q+1}^{1/2}M_L^{1/2}+\gamma_{q+1}^{1/2})\Big{(}\frac{M_0}2+C\lambda_{q+1}^{(6d+12)\alpha-\frac{d+1}{N}}+C(\delta_{q+1}^{1/2}M_L^{1/2}+\gamma_{q+1}^{1/2})\lambda_{q+1}^{(8d+8)\alpha-\frac{2d+1}{2N}}\Big{)}\notag\\
&\leq\frac34M_0(\delta_{q+1}^{1/2}M_L^{1/2}+\gamma_{q+1}^{1/2}),\label{2.19}
\end{align}
where we used conditions on the parameters to deduce
$$C\lambda_{q+1}^{(6d+12)\alpha-\frac{d+1}{N}}\leq\frac{M_0}8,\ C(\delta_{q+1}^{1/2}M_L^{1/2}+\gamma_{q+1}^{1/2})\lambda_{q+1}^{(8d+8)\alpha-\frac{2d+1}{2N}}\leq C\lambda_{q+1}^{(8d+9)\alpha-\frac{2d+1}{2N}} \leq\frac{M_0}8.$$
The above inequality together with (\ref{3.4}) yields for $t\in(4\sigma _q\wedge T_L,T_L]$
\begin{align*}
\|v_{q+1}(t)-v_q(t)\|_{L^2}&\leq\|v_{q+1}(t)-v_l(t)\|_{L^2}+\|v_{l}(t)-v_q(t)\|_{L^2}\notag\\
&\leq\|w_{q+1}(t)\|_{L^2}+l\|v_q(t)\|_{C_{t,x}^1}\notag\\
&\leq \frac34M_0(M_L^{1/2}\delta_{q+1}^{1/2}+\gamma_{q+1}^{1/2})+l\lambda_q^{d+1}M_L^{1/2}\\
&\leq M_0(M_L^{1/2}\delta_{q+1}^{1/2}+\gamma_{q+1}^{1/2}),
\end{align*}
where we used $-\frac{3\alpha}{2}+\frac{d}{b}<-\beta$ to deduce $l\lambda_q^{d+2}\leq\frac{M_0}{4}\delta_{q+1}^{1/2}$.

Similarly for $t\in(\frac{\sigma_q}{2}\wedge T_L,4\sigma _q\wedge T_L]$ we obtain
\begin{align}
\|w_{q+1}(t)\|_{L^2}&\leq \frac34M_0((M_L+qA)^{1/2}+\gamma_{q+1}^{1/2}),\label{4.25'}\\
\|v_{q+1}(t)-v_q(t)\|_{L^2}&\leq M_0((M_L+qA)^{1/2}+\gamma_{q+1}^{1/2}).
\end{align}

For $t\in[0,\frac{\sigma_q}{2}\wedge T_L]$, it holds $\chi(t)=0$ and $v_q(t)=0$, which implies that $\|v_{q+1}(t)-v_q(t)\|_{L^2}=0$. Hence (\ref{3.8}) follows and by a similar argument as \cite[Section 5.1.4]{HZZ} (\ref{3.3}) holds at the level $q+1$.

Next we estimate $C_{t,x}^i$-norm $i=1,2$, by (\ref{int2})-(\ref{int4}) and (\ref{3.13}) we have for $t\in[0,T_L]$
\begin{align}
\|w_{q+1}^{(p)}\|_{C_{t,x}^1}&\lesssim\sum_{\xi\in\Lambda} \|a_{(\xi)}\|_{C_{t,x}^{1}}\|W_{(\xi)}\|_{C_{t,x}^{1}}\notag\\
&\lesssim((M_L+qA)^{1/2}+\gamma_{q+1}^{1/2})l^{-3d-6}\lambda_{q+1}\mu\frac{r_\perp}{r_\parallel} r_\parallel^{-\frac{1}{2}}r_\perp^{-\frac{d-1}{2}}\notag\\
&\lesssim
((M_L+qA)^{1/2}+\gamma_{q+1}^{1/2})\lambda_{q+1}^{(6d+12)\alpha+1+d-\frac{4d+3}{2N}},\label{2.20}\\
\|w_{q+1}^{(p)}\|_{C_{t,x}^2}&\lesssim\sum_{\xi\in\Lambda} \|a_{(\xi)}\|_{C_{t,x}^{2}}\|W_{(\xi)}\|_{C_{t,x}^{2}}\notag\\
&\lesssim((M_L+qA)^{1/2}+\gamma_{q+1}^{1/2})l^{-4d-10}\lambda_{q+1}^2\mu^2(\frac{r_\perp}{r_\parallel})^2 r_\parallel^{-\frac{1}{2}}r_\perp^{-\frac{d-1}{2}}\notag\\
&\lesssim
((M_L+qA)^{1/2}+\gamma_{q+1}^{1/2})\lambda_{q+1}^{(8d+20)\alpha+\frac{3d}{2}+2-\frac{6d+5}{2N}},\label{2.21}
\end{align}

\begin{align}
\|w_{q+1}^{(c)}\|_{C_{t,x}^1}
&\lesssim\sum_{\xi\in\Lambda}\frac{1}{\lambda_{q+1}^2}\| a_{(\xi)}\|_{C_{t,x}^1}\|\nabla\Phi_{(\xi)}\xi\cdot\nabla\psi_{(\xi)}\|_{C_{t,x}^1}+\sum_{\xi\in\Lambda}\|\nabla a_{(\xi)}\|_{C_{t,x}^1}\|V_{(\xi)}\|_{C_{t,x}^1}\notag\\
&\lesssim ((M_L+qA)^{1/2}+\gamma_{q+1}^{1/2})l^{-4d-10}\lambda_{q+1}\mu r_\parallel^{-\frac{1}{2}}r_\perp^{-\frac{d-1}{2}}\notag\\
&\lesssim((M_L+qA)^{1/2}+\gamma_{q+1}^{1/2})\lambda_{q+1}^{(8d+20)\alpha+d+1-\frac{2d+1}{2N}},\label{2.21}
\end{align}
\begin{align}
\|w_{q+1}^{(c)}\|_{C_{t,x}^2}
&\lesssim\sum_{\xi\in\Lambda}\frac{1}{\lambda_{q+1}^2}\| a_{(\xi)}\|_{C_{t,x}^2}\|\nabla\Phi_{(\xi)}\xi\cdot\nabla\psi_{(\xi)}\|_{C_{t,x}^2}+\sum_{\xi\in\Lambda}\|\nabla a_{(\xi)}\|_{C_{t,x}^2}\|V_{(\xi)}\|_{C_{t,x}^2}\notag\\
&\lesssim ((M_L+qA)^{1/2}+\gamma_{q+1}^{1/2})l^{-5d-14}\lambda_{q+1}^2\mu^2r_\parallel^{-\frac{1}{2}}r_\perp^{-\frac{d-1}{2}}\notag\\
&\lesssim((M_L+qA)^{1/2}+\gamma_{q+1}^{1/2})\lambda_{q+1}^{(10d+28)\alpha+\frac{3d}{2}+2-\frac{2d+1}{2N}},\label{2.21}
\end{align}
and
\begin{align}
\|w_{q+1}^{(t)}\|_{C_{t,x}^1}&\lesssim\frac{1}{\mu}\sum_{\xi\in\Lambda}\left(\|a^2_{(\xi)}\psi^2_{(\xi)}\phi^2_{(\xi)}\|_{C_tW^{1+\alpha,\infty}}+\|a^2_{(\xi)}\psi^2_{(\xi)}\phi^2_{(\xi)}\|_{C_t^1W^{\alpha,\infty}}\right)\notag\\
&\lesssim(M_L+qA+\gamma_{q+1})l^{-5d-8}\lambda_{q+1}^{1+\alpha} r_\parallel^{-2}r_\perp^{-d+2}\notag\\
&\lesssim(M_L+qA+\gamma_{q+1})\lambda_{q+1}^{(10d+17)\alpha+1+d-\frac{3d+2}{N}}\label{2.22},
\end{align}
\begin{align}
\|w_{q+1}^{(t)}\|_{C_{t,x}^2}&\lesssim\frac{1}{\mu}\sum_{\xi\in\Lambda}\left(\|a^2_{(\xi)}\psi^2_{(\xi)}\phi^2_{(\xi)}\|_{C_tW^{2+\alpha,\infty}}+\|a^2_{(\xi)}\psi^2_{(\xi)}\phi^2_{(\xi)}\|_{C_t^1W^{1+\alpha,\infty}}+\|a^2_{(\xi)}\psi^2_{(\xi)}\phi^2_{(\xi)}\|_{C_t^2W^{\alpha,\infty}}\right)\notag\\
&\lesssim(M_L+qA+\gamma_{q+1})l^{-6d-12}\lambda_{q+1}^{2+\alpha}\mu r_\parallel^{-3}r_\perp^{-d+3}\notag\\
&\lesssim(M_L+qA+\gamma_{q+1})\lambda_{q+1}^{(12d+25)\alpha+\frac{3d}{2}+2-\frac{4d+3}{N}}.\label{2.26*}
\end{align}

Here we have a extra $\alpha$ since $\mathbb{P}\mathbb{P}_{\neq0}$ is not a bounded operator on $C^0$. In particular, we see that the fact that the time derivative of $\chi$ behaves like $2^{q+1}\lesssim l^{-1}$ does not pose any problems as the $C_{t,x}^0,C_{t,x}^1$-norms of $\tilde{w}_{q+1}^{(p)},\tilde{w}_{q+1}^{(c)},\tilde{w}_{q+1}^{(t)}$ always contain smaller powers of $l^{-1}$. Thus by $(12d+27)\alpha<\frac{1}{2N}$ and $M_L+qA+K\leq l^{-1}$ we obtain
$$\|v_{q+1}\|_{C_{t,x}^1}\leq\|v_{l}\|_{C_{t,x}^1}+\|w_{q+1}\|_{C_{t,x}^1}\leq M_L^{1/2}\lambda_{q+1}^{d+1},$$
$$\|v_{q+1}\|_{C_{t,x}^2}\leq\|v_{l}\|_{C_{t,x}^2}+\|w_{q+1}\|_{C_{t,x}^2}\leq M_L^{1/2}\lambda_{q+1}^{\frac{3d}{2}+2}.$$
This implies (\ref{3.4}) and (\ref{3.5}) at the level $q+1$.

We conclude this part with further $W^{1,p}$-norm for $p\in(1,\infty)$. By (\ref{int2})-(\ref{int4}) and (\ref{3.13})
\begin{align}
\|w_{q+1}^{(p)}(t)\|_{W^{1,p}}&\lesssim\sum_{\xi\in\Lambda}\|a_{(\xi)}\|_{C_{t,x}^{1}}\|W_{(\xi)}\|_{C_tW^{1,p}}\notag\\
&\lesssim ((M_L+qA)^{1/2}+\gamma_{q+1}^{1/2})l^{-3d-6} r_\perp^{\frac{d-1}{p}-\frac{d-1}{2}} r_\parallel^{\frac{1}{p}-\frac{1}{2}}\lambda_{q+1}\notag\\
&\lesssim((M_L+qA)^{1/2}+\gamma_{q+1}^{1/2})\lambda_{q+1}^{(6d+12)\alpha+(\frac{1}{p}-\frac{1}{2})(\frac{2d+1}{N}-d)+1},\label{2.27}
\end{align}
\begin{align}
\|w_{q+1}^{(c)}(t)\|_{W^{1,p}}&\lesssim\sum_{\xi\in\Lambda}\frac{1}{\lambda^2_{q+1}}\| a_{(\xi)}\|_{C_{t,x}^1}\|\nabla\Phi_{(\xi)}\xi\cdot\nabla\psi_{(\xi)}\|_{W^{1,p}}+\sum_{\xi\in\Lambda}\|a_{(\xi)}\|_{C_{t,x}^2}\|V_{(\xi)}\|_{W^{1,p}}\notag\\
&\lesssim((M_L+qA)^{1/2}+\gamma_{q+1}^{1/2})l^{-4d-10} \lambda_{q+1} r_\perp^{\frac{d-1}{p}-\frac{d-1}{2}} r_\parallel^{\frac{1}{p}-\frac{1}{2}}\notag\\
&\lesssim((M_L+qA)^{1/2}+\gamma_{q+1}^{1/2})\lambda_{q+1}^{(8d+20)\alpha+(\frac{1}{p}-\frac{1}{2})(\frac{2d+1}{N}-d)+1},\label{2.27*}
\end{align}
\begin{align}
\|w_{q+1}^{(t)}(t)\|_{W^{1,p}}&\lesssim\frac{1}{\mu}\sum_{\xi\in\Lambda}\|a_{(\xi)}^2\|_{C_{t,x}^{1}}\|\psi^2_{(\xi)}\phi^2_{(\xi)}\|_{C_tW^{1,p}}\notag\\
&\lesssim(M_L+qA+\gamma_{q+1})l^{-5d-8}\mu^{-1}r_\perp^{\frac{d-1}{p}-d+1}r_\parallel^{\frac{1}{p}-1}\lambda_{q+1}\notag\\
&\lesssim (M_L+qA+\gamma_{q+1})\lambda_{q+1}^{(10d+16)\alpha+1-\frac{d}{2}+(\frac{1}{p}-1)(\frac{2d+1}{N}-d)}.\label{2.28}
\end{align}
Now we recall $P(x)=(\frac{1}{x}-\frac{1}{2})(d-\frac{2d+1}{N})-1$ and we can choose a $\epsilon>0$ small enough such that $-P(s+\epsilon)\leq\alpha-P(s)$. Finally by $\|\chi(t)\|_{L^\infty}\leq1$ we obtain
\begin{align}\label{ww1s}
\|{w}_{q+1}\|_{C_tW^{1,s}}&\lesssim\|w_{q+1}\|_{C_tW^{1,s+\epsilon}}\lesssim \lambda_{q+1}^{(10d+16)\alpha-P(s+\epsilon)}\notag\\
&\lesssim\lambda_{q+1}^{(10d+17)\alpha-P(s)} \leq \frac {1}2M_L^{1/2}\delta_{q+1}^{1/2}.
\end{align}
Here we used the conditions on the parameters in the second and the last inequalities.
Hence we obtain
\begin{align*}
\|v_{q+1}(t)-v_q(t)\|_{W^{1,s}}&\leq\|{w}_{q+1}(t)\|_{W^{1,s}}+l\|v_q\|_{C_{t,x}^2}\\
&\leq \frac12M_L^{1/2}\delta_{q+1}^{1/2}+l\lambda_q^{\frac{3d}{2}+2}M_L^{1/2}\leq M_L^{1/2}\delta_{q+1}^{1/2},
\end{align*}
which implies (\ref{3.10}). Here we used our conditions on parameters to deduce $\beta+\frac{3d}{2b}< \frac{3\alpha}{2}$ and $l\lambda_q^{\frac{3d}{2}+2}\leq\frac12\delta_{q+1}^{1/2}$.

\subsubsection{Proof of (\ref{3.11})}\label{subsub5}
 For $t\in(4\sigma_q\wedge T_L,T_L]$ we have $\chi(t)=1$, then
\begin{align}
|\|v_{q+1}\|_{L^2}^2-\|v_q\|_{L^2}^2-d\gamma_{q+1}|&\leq| \|w_{q+1}^{(p)}\|_{L^2}^2-d\gamma_{q+1}|+\|w_{q+1}^{(c)}+w_{q+1}^{(t)}\|_{L^2}^2
+2\|v_l(w_{q+1}^{(c)}+w_{q+1}^{(t)})\|_{L^1}\notag\\
&+2\|v_lw_{q+1}^{(p)}\|_{L^1}+2\|w_{q+1}^{(p)}(w_{q+1}^{(c)}+w_{q+1}^{(t)})\|_{L^1}+|\|v_l\|_{L^2}^2-\|v_q\|_{L^2}^2|.\label{4.34}
\end{align}
We use (\ref{2.12}) and the fact that $\mathring{R}_l$ is traceless to deduce
\begin{align*}
|w_{q+1}^{(p)}|^2-\frac{d\gamma_{q+1}}{(2\pi)^d}&=2d\sqrt{l^2+|\mathring{R}_l|^2}
+\sum_{\xi\in\Lambda}a_{(\xi)}^2\mathbb{P}_{\neq0}|W_{(\xi)}|^2,
\end{align*}
hence
\begin{align}
| \|w_{q+1}^{(p)}\|_{L^2}^2-d\gamma_{q+1}|&\leq 2d\cdot(2\pi)^dl+2d\|\mathring{R}_l\|_{L^1}+\sum_{\xi\in\Lambda}|\int_{\mathbb{T}^d}a_{(\xi)}^2\mathbb{P}_{\neq0}|W_{(\xi)}|^2\dif x|\notag.
\end{align}
Now we estimate each term separately. Using (\ref{3.6}) and $\supp \phi_l\subset[0,l]$
$$2d\|\mathring{R}_l(t)\|_{L^1}\leq2d\delta_{q+1}M_L,$$
and
$$2d\cdot(2\pi)^dl\leq2d\cdot(2\pi)^d\lambda_{q+1}^{-\frac{3\alpha}{2}}\leq\frac{1}{7}\delta_{q+1}M_L,$$
which requires $2\beta<\frac{3\alpha}{2}$ and we choose $a$ large enough to absorb the constant.

We obverse that $W_{(\xi)}$ is $(\mathbb{T}/r_\perp\lambda_{q+1})^d$-periodic so
$$\mathbb{P}_{\neq0}(W_{(\xi)}\otimes W_{(\xi)})=\mathbb{P}_{\geq\frac{r_\perp\lambda_{q+1}}{2}}(W_{(\xi)}\otimes W_{(\xi)}),$$
where $\mathbb{P}_{\geq r}=\Id-\mathbb{P}_{< r}$, and $\mathbb{P}_{<r}$ denotes the Fourier multiplier operator, which projects a function
onto its Fourier frequencies $< r$ in absolute value. By \cite[Proposition C.1]{TCP}  we obtain
\begin{align}
\sum_{\xi\in\Lambda}|\int_{\mathbb{T}^d}a_{(\xi)}^2\mathbb{P}_{\neq0}|W_{(\xi)}|^2
\dif x|&\lesssim\sum_{\xi\in\Lambda}|\int_{\mathbb{T}^d}a_{(\xi)}^2\mathbb{P}_{\geq\frac{r_\perp\lambda_{q+1}}{2}}|W_{(\xi)}|^2
\dif x|\notag\\
&\lesssim\sum_{\xi\in\Lambda}|\int_{\mathbb{T}^d}|\nabla|^{N^*}a_{(\xi)}^2|\nabla|^{-N^*}\mathbb{P}_{\geq\frac{r_\perp\lambda_{q+1}}{2}}|W_{(\xi)}|^2
\dif x|\notag\\
&\lesssim\|a_{(\xi)}^2\|_{C^{N^*}}(r_\perp\lambda_{q+1})^{-N*}\|\ |W_{(\xi)}|^2\|_{L^2}\notag\\
&\lesssim\lambda_{q+1}^{((2d+8)N^*+(8d+8))\alpha-\frac{N^*}{N}-\frac{2d+1}{2N}+\frac{d}{2}}(M_L+K)\notag\\
&\lesssim\lambda_{q+1}^{((2d+8)N^*+(8d+8))\alpha-\frac{N^*}{N}-\frac{2d+1}{2N}+\frac{d}{2}}\lambda_{q+1}^{2\alpha}\leq\frac{1}{7}\delta_{q+1}M_L,\label{N*}
\end{align}
where we used $M_L+K\leq l^{-1}\leq\lambda_{q+1}^{2\alpha}$. Here we may choose $N^*=[\frac{dN}2]+2$ with $N$ given in Section \ref{Cp} and then $-\frac{N^*}{N}+\frac{d}{2}<-\frac{1}{N}<0$. The last inequality holds since by our choice of parameters we have
$$(2d+8)N^*\alpha <\frac{1}{N},\ \ (8d+10)\alpha+2\beta<\frac{2d+1}{2N}.$$

Go back to (\ref{4.34}), by (\ref{2.18}) and (\ref{2.18*})
\begin{align}
\|w_{q+1}^{(c)}+w_{q+1}^{(t)}\|_{L^2}^2
&\lesssim (M_L+\gamma_{q+1})(\lambda_{q+1}^{(6d+12)\alpha-\frac{d+1}N}
+\lambda_{q+1}^{(8d+9)\alpha-\frac{2d+1}{2N}})^2\leq\frac{1}{7}\delta_{q+1}M_L,\notag
\end{align}
where we used $M_L+K\leq l^{-1}$ and the conditions on the parameters to deduce $(16d+20)\alpha-\frac{2d+1}{N}<-2\beta b$ in the last inequality.

By (\ref{3.3}), (\ref{2.16''}), (\ref{2.18}) and (\ref{2.18*})
\begin{align}
2\|v_l&(w_{q+1}^{(c)}+w_{q+1}^{(t)})\|_{L^1}
+2\|w_{q+1}^{(p)}(w_{q+1}^{(c)}+w_{q+1}^{(t)})\|_{L^1}\notag\\
&\lesssim M_0((M_L+qA)^{1/2}+K^{1/2})\|w_{q+1}^{(c)}+w_{q+1}^{(t)}\|_{L^2}\notag\\
&\lesssim M_0((M_L+qA)^{1/2}+K^{1/2})(M_L^{1/2}+\gamma_{q+1}^{1/2})(\lambda_{q+1}^{(6d+12)\alpha-\frac{d+1}N}
+\lambda_{q+1}^{(8d+9)\alpha-\frac{2d+1}{2N}})\notag\\
&\lesssim M_0(\lambda_{q+1}^{(6d+14)\alpha-\frac{d+1}N}
+\lambda_{q+1}^{(8d+11)\alpha-\frac{2d+1}{2N}})
\leq\frac{2}{7}\delta_{q+1}M_L,\notag
\end{align}
where we used $M_L+qA+K\leq l^{-1}$ in the last second inequality and in the last inequality we used the conditions on the parameters to deduce $(8d+11)\alpha-\frac{2d+1}{2N}<-2\beta b$ and $a$ large enough to absorb the constant $M_0$.

By (\ref{3.4}) and (\ref{2.17})
\begin{align}
2\|v_lw_{q+1}^{(p)}\|_{L^1}&\lesssim\|v_l\|_{L^\infty}\|w_{q+1}^{(p)}\|_{L^1}\notag\\
&\lesssim\lambda_q^{d+1}M_L^{1/2}(M_L^{1/2}\delta_{q+1}^{1/2}+\gamma_{q+1}^{1/2})\lambda_{q+1}^{(4d+4)\alpha+\frac{1}{2}(\frac{2d+1}{N}-d)}\notag\\
&\lesssim\lambda_{q+1}^{(4d+7)\alpha+\frac{1}{2}(\frac{2d+1}{N}-d)}\leq\frac{1}{7}\delta_{q+1}M_L,\notag
\end{align}
where $M_L+qA+K\leq l^{-1}$ and the consitions on the parameters to deduce $\lambda_q^{d+1}\leq\lambda_{q+1}^\alpha
$ and $(4d+7)\alpha+\frac{1}{2}(\frac{2d+1}{N}-d)<-2\beta$.

By (\ref{3.3}) and (\ref{2.17})
\begin{align}
|\|v_l\|_{L^2}^2-\|v_q\|_{L^2}^2|&\leq\|v_l-v_q\|_{L^2}(\|v_l\|_{L^2}+\|v_q\|_{L^2})\notag\\
&\lesssim l\lambda_{q}^{d+1}M_L^{1/2}M_0(M_L+qA+K)^{1/2}\notag\\
&\leq\frac1{7}\lambda_{q+1}^{-2\beta}M_L\leq\frac17\delta_{q+1}M_L,\notag
\end{align}
where we used the conditions on the parameters to deduce $\lambda_q^{d+1}\leq\lambda_{q+1}^{\alpha/2}$, $M_L+qA+K\leq\lambda_{q+1}^{\alpha-2\beta}$ and $a$ large enough to absorb the extra constant.

Combining the above estimates (\ref{3.11}) follows.

\subsubsection{Construction of the Reynolds stress}\label{cotrs}
From (\ref{3.1}) and (\ref{eq:v_l}) we obtain

\begin{align}
\div&\mathring{R}_{q+1}-\nabla\pi_{q+1}+\nabla\pi_{l}\notag\\
&=\partial_t(\tilde{w}_{q+1}^{(p)}+\tilde{w}_{q+1}^{(c)})+\div((v_l+z_l)\otimes w_{q+1})+\div(w_{q+1}\otimes (v_l+z_l))\ (:=\div R_{lin}+\nabla \pi_{lin})\notag\\
&+\div((\tilde{w}_{q+1}^{(c)}+\tilde{w}_{q+1}^{(t)})\otimes w_{q+1}+\tilde{w}_{q+1}^{(p)}\otimes ((\tilde{w}_{q+1}^{(c)}+\tilde{w}_{q+1}^{(t)}))\ (:=\div R_{cor}+\nabla \pi_{cor})\notag\\
&+\partial_t\tilde{w}_{q+1}^{(t)}+\div(\tilde{w}_{q+1}^{(p)}\otimes \tilde{w}_{q+1}^{(p)}+\mathring{R}_l)\ (:=\div R_{osc}+\nabla \pi_{osc})\notag\\
&-\div(\mathcal{A}(Dv_l+Dz_{q+1}+Dw_{q+1})-\mathcal{A}(Dv_l+Dz_l))-\Delta z_l+\Delta z_{q+1}\ (:=\div{R}_{nonlin}+\nabla \pi_{nonlin})\notag\\
&+\div(R_{com1}+{R}_{com2})\notag\\
&+\div((z_{q+1}-z_l)\otimes v_{q+1}+v_{q+1}\otimes (z_{q+1}-z_l)+z_{q+1}\otimes z_{q+1}-z_l\otimes z_l)(:=\div R_{com3}+\nabla\pi_{com3}),
\end{align}
where
\begin{align}
R_{com1}&=(v_l+z_l)\mathring{\otimes}(v_l+z_l)-((v_q+z_q)\mathring{\otimes}(v_q+z_q))*_x\phi_l*\varphi_l,\notag \\
R_{com2}&=(\mathcal{A}(Dv_q+Dz_q)*_x\phi_l)*\varphi_l-\mathcal{A}(Dv_l+Dz_l).\notag
\end{align}

By using $\mathcal{R}$ introduced in Section \ref{tamr} we define
\begin{align*}
R_{lin}:&=\mathcal{R}(\partial_t(\tilde{w}_{q+1}^{(p)}+\tilde{w}_{q+1}^{(c)}))+(v_l+z_l)\mathring{\otimes} w_{q+1}+w_{q+1}\mathring{\otimes} (v_l+z_l),    \\
R_{cor}:&=(\tilde{w}_{q+1}^{(c)}+\tilde{w}_{q+1}^{(t)})\mathring{\otimes} w_{q+1}+\tilde{w}_{q+1}^{(p)}\mathring{\otimes} (\tilde{w}_{q+1}^{(c)}+\tilde{w}_{q+1}^{(t)}),\\
R_{com3}:&=(z_{q+1}-z_l)\mathring{\otimes} v_{q+1}+v_{q+1}\mathring{\otimes} (z_{q+1}-z_l)+z_{q+1}\mathring{\otimes} z_{q+1}-z_l\mathring{\otimes} z_l,\\
R_{nonlin}:&=-\mathcal{A}(Dv_l+Dz_{q+1}+Dw_{q+1})+\mathcal{A}(Dv_l+Dz_l)-2D z_l+2D z_{q+1}.
\end{align*}
Here similar to $R_{com2}$, the nonlinear error $R_{nonlin}$ is trace-free.

In order to define the remaining oscillation error in the forth line, similarly as \cite[Section 7.6.1]{BV}, we apply (\ref{2.12}) and (\ref{int1}) to obtain
\begin{align*}
\div R_{osc}&=\partial_t\tilde{w}_{q+1}^{(t)}+\div(\tilde{w}_{q+1}^{(p)}\otimes \tilde{w}_{q+1}^{(p)}+\mathring{R}_l)\\
=&-\sum_{\xi\in\Lambda}\chi^2 \mathbb{P}_{\neq0}\left(\nabla a_{(\xi)}^2\mathbb{P}_{\neq0} (W_{(\xi)}\otimes W_{(\xi)})\right)-\chi^2\frac{1}{\mu}\sum_{\xi\in\Lambda}\mathbb{P}_{\neq0}\left(\partial_ta_{(\xi)}^2(\phi_{(\xi)}^2\psi_{(\xi)}^2\xi)\right)+\chi^2\nabla p+\chi^2\nabla \rho\\
&+(\chi^2)'w_{q+1}^{(t)}+\div( (1-\chi^2)\mathring{R}_{l}).
\end{align*}
Therefore
\begin{align*}
R_{osc}:&=\sum_{\xi\in\Lambda}\chi^2\mathcal{B}\left(\nabla a_{(\xi)}^2,\mathbb{P}_{\neq0}(W_{(\xi)}\otimes W_{(\xi)})\right)-\frac{1}{\mu}\sum_{\xi\in\Lambda}\chi^2\mathcal{R}\left((\partial_ta_{(\xi)}^2)\phi_{(\xi)}^2\psi_{(\xi)}^2\xi\right)+(\chi^2)'\mathcal{R} w_{q+1}^{(t)}+ (1-\chi^2)\mathring{R}_{l}\\
:&=R_{osc}^{(x)}+R_{osc}^{(t)}+(\chi^2)'\mathcal{R}w_{q+1}^{(t)}+(1-\chi^2)\mathring{R}_{l},
\end{align*}
where $\mathcal{B}$ is introduced in Section \ref{tamr}.

Finally define the Reynolds stress on the level $q+1$ by $$\mathring{R}_{q+1}=R_{lin}+R_{cor}+R_{osc}+
{R}_{nonlin}+R_{com1}+{R}_{com2}+R_{com3}.$$
It is easy to see that $\mathring{R}_{q+1}$ is a trace-free and symmetric matrix.
\subsubsection{Estimate of $\mathring{R}_{q+1}$}\label{subsub7}
To conclude the proof of Proposition \ref{prop:1} we shall verify (\ref{3.9}).
In order to establish the iterative estimate, we distinguish three cases corresponding to the three time intervals.

\textbf{1}. Let $t\in(\sigma_q\wedge T_L,T_L]$.
Note that if $T_L\leq \sigma_q$ then there is nothing to estimate here, hence
we assume that $\sigma_q < T_L$ and $t\in(\sigma_q, T_L]$. Then we have $\chi(t)=1$. We estimate each term in the definition of $R_{q+1}$ separately.

For the linear error $R_{lin}$, by (\ref{int4}), (\ref{***}) and (\ref{3.13}) we obtain for $\epsilon>0$ small enough

\begin{align}\label{rlin1}
\|\mathcal{R}(\partial_t(w_{q+1}^{(p)}+w_{q+1}^{(c)}))(t)\|_{L^1}&\lesssim\sum_{\xi\in\Lambda}\|\partial_t (a_{(\xi)}V_{(\xi)})\|_{C_tL^{1+\epsilon}}\\\notag
&\lesssim ((M_L+qA)^{1/2}+\gamma_{q+1}^{1/2})l^{-3d-6}r_\perp^{\frac{d-1}{1+\epsilon}-\frac{d-1}{2}} r_\parallel^{\frac{1}{1+\epsilon}-\frac{1}{2}}\frac{r_\perp\mu}{r_\parallel}\\\notag
&\lesssim\lambda_{q+1}^{(6d+13)\alpha-\frac{1}{2N}+\frac{\epsilon}{1+\epsilon}(d-\frac{2d+1}{N})}\lesssim\lambda_{q+1}^{(6d+14)\alpha-\frac{1}{2N}},
\end{align}
where we choose $\epsilon>0$ small enough such that $d \epsilon<\alpha$.

By \cite[Lemma 9]{DV}
$$\|z^{in}(t)\|_{L^\infty}=\|e^{t\Delta}u_0\|_{L^\infty}\lesssim(1+t^{-d/4})\|u_0\|_{L^2}. $$
Hence
by (\ref{3.2}), (\ref{2.17})-(\ref{2.18*})
\begin{align*}
\|w_{q+1}\mathring{\otimes}z_l(t)\|_{L^1}&\leq \|z_l(t)\|_{L^\infty}\|w_{q+1}(t)\|_{L^1}\lesssim(\sup_{s\in[t-l,t]}s^{-d/4}+\lambda_{q+1}^{\alpha/4})M_L^{1/2}\|w_{q+1}(t)\|_{L^{1+\epsilon}}\\
&\lesssim(l^{-d/4}+\lambda_{q+1}^{\alpha/4})M_L^{1/2}((M_L+qA)^{1/2}+\gamma_{q+1}^{1/2})\lambda_{q+1}^{(4d+4)\alpha-\frac{d}{2}+\frac{2d+1}{2N}+\epsilon d}\\
&\lesssim \lambda_{q+1}^{(\frac{9}{2}d+7)\alpha-\frac{d}{2}+\frac{2d+1}{2N}}\lesssim\lambda_{q+1}^{(\frac{9}{2}d+7)\alpha-\frac{d}{4}},
\end{align*}
where we choose $\epsilon>0$ small enough such that $ d \epsilon<\alpha$. Here we used $t>\sigma_q\geq 2l$, $N>4$ and $M_L+qA+K\leq l^{-1}$. Thus by our choice of parameters we obtain
\begin{align}
\|R_{lin}(t)\|_{L^1}\leq \lambda_{q+1}^{(6d+14)\alpha-\frac{1}{2N}}+\lambda_{q+1}^{(\frac{9}{2}d+7)\alpha-\frac{d}4}\leq \frac{M_L\delta_{q+2}}{7}.\label{3.15}
\end{align}

The corrector error is estimated using (\ref{2.18+-})-(\ref{2.18**}) and (\ref{4.25'}) as
\begin{align}
\|R_{cor}(t)\|_{L^1}&\leq\|w_{q+1}^{(c)}(t)+w_{q+1}^{(t)}(t)\|_{L^{2}}\|w_{q+1}(t)\|_{L^{2}}
+\|w_{q+1}^{(c)}(t)+w_{q+1}^{(t)}(t)\|_{L^{2}}\|w_{q+1}^{(p)}(t)\|_{L^{2}}\notag\\
&\lesssim((M_L+qA)^{1/2}+\gamma_{q+1}^{1/2})(\lambda_{q+1}^{(6d+12)\alpha-\frac{d+1}N}
+\lambda_{q+1}^{(8d+8)\alpha-\frac{2d+1}{2N}})\times
\frac34M_0((M_L+qA)^{1/2}+K^{1/2})\notag\\
&\lesssim M_0\lambda_{q+1}^{(8d+10)\alpha-\frac{2d+1}{2N}}
\leq \frac{M_L\delta_{q+2}}{7},\label{2.31}
\end{align}
where we used $M_L+qA+k\leq l^{-1}$ and the conditions on the parameters to deduce $(8d+10)\alpha-\frac{2d+1}{2N}<-2\beta b$.

We continue with $R_{osc}$. For $t\in(\sigma_q\wedge T_L,T_L]$ we have $R_{osc}=R_{osc}^{(x)}+R_{osc}^{(t)}$. In order to bound the first term, we apply Theorem \ref{bb} and by (\ref{int4}) and (\ref{3.13})
\begin{align}
\|R_{osc}^{(x)}(t)\|_{L^1}&\lesssim\sum_{\xi\in\Lambda}\|\mathcal{B}(\nabla a_{(\xi)}^2,\mathbb{P}_{\neq0}(W_{(\xi)}\otimes W_{(\xi)}))\|_{C_tL^1}\lesssim \sum_{\xi\in\Lambda}\|\nabla a_{(\xi)}^2\|_{C_tC^1}\|\mathcal{R}(W_{(\xi)}\otimes W_{(\xi)})\|_{C_tL^{1+\epsilon}}\notag\\
&\lesssim l^{-6d-12}(M_L+qA+\gamma_{q+1})(r_\perp\lambda_{q+1})^{-1}\|W_{(\xi)}\otimes W_{(\xi)}\|_{C_tL^{1+\epsilon}}\notag\\
&\lesssim l^{-6d-12}(M_L+qA+\gamma_{q+1})(r_\perp\lambda_{q+1})^{-1}r_\perp^{(d-1)(\frac{1}{1+\epsilon}-1)}r_\parallel^{\frac{1}{1+\epsilon}-1}\notag\\
&\lesssim \lambda_{q+1}^{(12d+26)\alpha-\frac{1}{N}+d\epsilon}\lesssim\lambda_{q+1}^{(12d+27)\alpha-\frac{1}{N}}
\leq\frac{M_L\delta_{q+2}}{14},\label{2.32}
\end{align}
where we used $M_L+qA+K\leq l^{-1}$ and the conditions on the parameters to deduce $(12d+27)\alpha-\frac{1}{N}<-2\beta b$. Here we chose $\epsilon>0$ small enough such that $d\epsilon<\alpha$.

For the second term ${R}_{osc}^{(t)}$ we use (\ref{int2})-(\ref{int3}) and (\ref{3.13}) to deduce
\begin{align}
\|R_{osc}^{(t)}(t)\|_{L^1}&\lesssim \frac{1}{\mu}\|(\partial_ta_{(\xi)}^2)\phi_{(\xi)}^2\psi_{(\xi)}^2\xi\|_{C_tL^{1+\epsilon}}\lesssim(M_L+qA+\gamma_{q+1}) l^{-5d-8}\mu^{-1}r_\perp^{(d-1)(\frac{1}{1+\epsilon}-1)}r_\parallel^{\frac{1}{1+\epsilon}-1}\notag\\
&\lesssim \lambda_{q+1}^{(10d+18)\alpha-\frac{d}{2}+d\epsilon}
\leq\frac{M_L\delta_{q+2}}{14},\label{2.32}
\end{align}
where we used $M_L+qA+K\leq l^{-1}$ and  the conditions on the parameters to deduce $(10d+19)\alpha-\frac{d}{2}<-2\beta b$. Here we  chose $\epsilon>0$ small enough such that $d\epsilon<\alpha$

The nonlinear error will be divided into three cases. For $r\in(1,2],\nu_0>0$, by lemma \ref{lem:1}
\begin{align}
\|{R}_{nonlin}(t)\|_{L^1}
&\lesssim \|\nabla w_{q+1}(t)\|_{L^{1}}+\|\nabla z_{q+1}(t)-\nabla z_l(t)\|_{L^1}.
\end{align}
By (\ref{ww1s}) we have
\begin{align}
\|\nabla w_{q+1}(t)\|_{L^1}\lesssim\lambda_{q+1}^{(10d+17)\alpha-P(1)}
\leq\frac{M_L\delta_{q+2}}{14},\notag
\end{align}
where we used the conditions on the parameters to deduce $(10d+17)\alpha-P(1)<-2\beta b$. Moreover we have

\begin{align}
\|\nabla z_{q+1}(t)-\nabla z_l(t)\|_{L^{r^*}}&\lesssim\|\nabla z_{q+1}(t)-\nabla z_q(t)\|_{L^{r^*}}+\|\nabla z_q(t)-\nabla z_l(t)\|_{L^{r^*}}.\notag
\end{align}
We first estimate the second term,
\begin{align}
\|\nabla z_q(t)-\nabla z_l(t)\|_{L^{r^*}}
&\lesssim\|\nabla z^{in}(t)-\nabla z^{in}_l(t)\|_{L^{r^*}}+\|\nabla Z_q(t)-\nabla Z_l(t)\|_{L^2}\notag\\
&\lesssim l^{1/2}(\|\nabla z^{in}\|_{C_{[\frac{\sigma_q}{2},t]}H^{1/2}}+\|\nabla z^{in}\|_{C^{1/2}_{[\frac{\sigma_q}{2},t]}L^{r^*}})\notag\\
&\ \ \ \ \ +l\|\nabla^2Z_q\|_{C_{[\frac{\sigma_q}{2},t]}L^2}+l^{1/2-2\delta}\|\nabla Z_q\|_{C_{[\frac{\sigma_q}{2},t]}^{1/2-2\delta}L^2},\notag
\end{align}
where we denote $z^{in}_l:=(z^{in}*_x\phi_l)*_t\varphi_l,Z_l:=(Z_q*_x\phi_l)*_t\varphi_l$. By \cite[Lemma 9]{DV} we have
\begin{align}
\|\nabla z^{in}(t)\|_{H^{1/2}}\lesssim (1+t^{-3/4})\|u_0\|_{L^2},\label{zin1}
\end{align}
for $|t-s|\leq1$
\begin{align}
\|\nabla z^{in}(t)-\nabla z^{in}(s)\|_{L^{r^*}}&\lesssim\|\nabla e^{s\Delta}(e^{(t-s)\Delta}u_0-u_0)\|_{L^2}=\|\nabla e^{s\Delta}\int_0^{t-s}\Delta e^{r\Delta}u_0dr\|_{L^2}\notag\\
&\lesssim \int_0^{t-s}(1+r^{-1/2})\|(-\Delta)e^{s\Delta}u_0\|_{L^2}dr\lesssim |t-s|^{1/2}(1+s^{-1})\|u_0\|_{L^2},\notag
\end{align}
for $|t-s|>1$
\begin{align}
\|\nabla z^{in}(t)-\nabla z^{in}(s)\|_{L^{r^*}}&\leq2\|u_0\|_{W^{1,r^*}}\leq 2|t-s|^{1/2}\|u_0\|_{W^{1,s}}.\notag
\end{align}
Hence, we obtain
\begin{align}
\|\nabla z^{in}\|_{C^{1/2}_{[\frac{\sigma_q}{2},t]}L^{r^*}}\lesssim (1+\sigma_{q}^{-1})M_L^{1/2}.\label{zin4}
\end{align}
Thus by (\ref{3.2}), (\ref{zin1})-(\ref{zin4}) and our assumption $\sigma_q\geq l^{1/4}$ we obtain
\begin{align}
\|\nabla z_{q}(t)-\nabla z_l(t)\|_{L^{r^*}}&\lesssim l^{1/2}M_L^{1/2}(1+\sigma_q^{-3/4}+\sigma_q^{-1})+l^{1/2-2\delta}L\lambda_{q+1}^{\frac\alpha4}\notag\\
&\lesssim M_L^{1/2}(l^{1/4}+l^{1/2-2\delta}\lambda_{q+1}^{\frac\alpha4})\lesssim M_L^{1/2}\lambda_{q+1}^{-\alpha/8},\label{bdd:Dz1}
\end{align}
where we used the conditions on the parameters to deduce $l^{1/4}\lambda_{q+1}^{\frac{\alpha}8}+l^{1/2-2\delta}\lambda_{q+1}^{\frac{3\alpha}8}\leq1$.

Thus by (\ref{3.2}) and above estimate
\begin{align}
\|\nabla z_{q+1}(t)-\nabla z_l(t)\|_{L^{r^*}}
&\lesssim M_L^{1/2}(\lambda_{q+1}^{-\frac{\alpha\sigma}{2d+4}}+\lambda_{q+1}^{-\alpha/8}) \leq\frac{M_L\delta_{q+2}}{14},\label{bdd:Dz}
\end{align}
where we used our conditions on the parameters to deduce $\frac{\alpha}{8}>2\beta b,\frac{\alpha\sigma}{2d+4}>2\beta b$.

Similarly for $r\in(1,2],\nu_0=0$, we obtain
\begin{align}
\|{R}_{nonlin}(t)\|_{L^1}&\lesssim \|\nabla w_{q+1}(t)\|_{L^1}^{r-1}+\|\nabla z_{q+1}(t)-\nabla z_l(t)\|_{L^1}^{r-1}\notag\\
&\lesssim \lambda_{q+1}^{((10d+17)\alpha-P(1))(r-1)}+M_L^{(r-1)/2}(\lambda_{q+1}^{-\frac{\alpha\sigma}{2d+4}}+\lambda_{q+1}^{-\alpha/8})^{r-1} \leq\frac{M_L\delta_{q+2}}{7},
\end{align}
where we used the conditions on the parameters to deduce $(P(1)-(10d+17)\alpha)(r-1)>2\beta b, \frac{\alpha}{8}(r-1)>2\beta b,\frac{\alpha\sigma}{2d+4}(r-1)>2\beta b$ and $a$ large enough to absorb $M_L$.

For the case $ r\in(2,\frac{3d+2}{d+2})$ by Lemma \ref{lem:1}, (\ref{ww1s}), (\ref{bdd:Dz}) and H\"older's inequality we have
\begin{align}
\|\mathring{R}_{nonlin}(t)\|_{L^1}&\lesssim(\|\nabla w_{q+1}(t)\|_{L^{r-1}}+\|\nabla z_{q+1}(t)-\nabla z_l(t)\|_{L^{r-1}})\notag\\
&(1+\|\nabla v_l(t)+\nabla z_{q+1}(t)+\nabla w_{q+1}(t)\|_{L^{r-1}}^{r-2}+\|\nabla v_l(t)+\nabla z_l(t)\|_{L^{r-1}}^{r-2}).\notag\\
&\lesssim(\lambda_{q+1}^{(10d+17)\alpha-P(r^*)}+M_L^{1/2}\lambda_{q+1}^{-\frac{\alpha\sigma}{2d+4}}+M_L^{1/2}\lambda_{q+1}^{-\alpha/8}) M_L^{1/2}(1+\lambda_q^{(d+1)(r-2)})\notag\\
&\lesssim \lambda_{q+1}^{(10d+18)\alpha-P(r^*)}\lambda_q^{d+1}+M_L(\lambda_{q+1}^{-\frac{\alpha\sigma}{2d+4}}+\lambda_{q+1}^{-\alpha/8})\lambda_q^{d+1}\notag\\
&\lesssim \lambda_{q+1}^{(10d+19)\alpha-P(r^*)}+M_L\lambda_{q+1}^{-\frac{\alpha}{16}}+M_L\lambda_{q+1}^{-\frac{\alpha}{2(2d+4)}\sigma}\leq\frac{M_L\delta_{q+2}}{7},\label{3.18}
\end{align}
where we used
$$\|\nabla v_q\|_{C_tL^{r-1}}^{r-2}\lesssim M_L^{1/2}
\lambda_q^{(d+1)(r-2)},\ \ \|\nabla z_{q+1}\|_{C_tL^{r-1}}^{r-2}\lesssim\|\nabla z_{q+1}\|_{C_tL^2}^{r-2}\leq M_L^{1/2}$$
and the conditions on the parameters to deduce $\alpha b>16(d+1)$, ${\alpha}>32\beta b$, $\frac{\alpha\sigma}{2(2d+4)}>2\beta b$, $P(r^*)-(10d+19)\alpha>2\beta b$, $\alpha b>\frac{2(d+1)(2d+4)}{\sigma}$ and $a$ large enough to absorb $M_L$.

Thus we obtain
\begin{align}
\|{R}_{nonlin}(t)\|_{L^1}\leq\frac{M_L\delta_{q+2}}{7}.\notag
\end{align}

We continue with $R_{com1}$. Using a standard mollification estimate we obtain
\begin{align}
\|R_{com1}(t)\|_{L^1}&\lesssim
(l\|v_q\|_{C_{x,t}^1}+l^{1/2-2\delta}\|Z_q\|_{C_t^{1/2-2\delta}L^2}+l^{1/2-2\delta}\|Z_q\|_{C_tH^{1/2-2\delta}})(\|v_q\|_{L_t^\infty L^2}+\|z_q\|_{C_tL^2})\notag\\
&+l^{1/2}(\|z^{in}\|_{C_{[\frac{\sigma_q}{2},t]}^{1/2}L^2}+\|z^{in}\|_{C_{[\frac{\sigma_q}{2},t]}H^{1/2}})(\|v_q\|_{L_t^\infty L^2}+\|z_q\|_{C_tL^2}).\notag
\end{align}
By  \cite[Lemma 9]{DV} we have
\begin{align}
\|z^{i n}(t)\|_{H^{1 / 2}} \lesssim(t^{-1 / 4}+1)\|u_{0}\|_{L^{2}},\label{zin11}
\end{align}
and similarly as (\ref{zin4})
\begin{align}
\|z^{in}\|_{C_{[\frac{\sigma_q}{2},t]}^{1/2}L^2}\leq (1+\sigma_q^{-1/2})M_L^{1/2}.\label{zin13}
\end{align}
Hence, together with (\ref{3.2})-(\ref{3.4}) we obtain
\begin{align}
\|R_{com1}(t)\|_{L^1}&\lesssim(l\lambda_q^{d+1}+l^{1/2-2\delta}+l^{1/2}(1+\sigma_q^{-1/2}))M_L^{\frac12}\times M_0(M_L+qA+K)^{1/2}\notag\\
&\lesssim \lambda_{q+1}^{-\alpha/8}M_0(M_L+qA+K)^{1/2}M_L^{\frac12}\leq\frac{M_L\delta_{q+2}}{7},\label{3.19}
\end{align}
where we used the conditions on the parameters to deduce $l\lambda_{q}^{d+1}\leq \lambda_{q+1}^{-\frac{\alpha}{8}}$, $M_0(M_L+qA+K)\leq\frac{M_L}{7} \lambda_{q+1}^{\alpha/8-2\beta b}$.

We write ${R}_{com2}$ as follows:
\begin{align}
R_{com2}&=\mathcal{A}(Dv_q+Dz_q)-(\mathcal{A}(Dv_q+Dz_q)*_x\phi_l)*_t\varphi_l\ \ (:={R}_{com21})\notag\\
&+\mathcal{A}(Dv_l+Dz_l)-\mathcal{A}(Dv_q+Dz_q)\ \ (:={R}_{com22}).\notag
\end{align}
First we deal with ${R}_{com21}$. For $r\in(1,2],\nu_0>0$, by lemma \ref{lem:1}
\begin{align}
\|{R}_{com21}(t)\|_{L^1}
&\lesssim l\|\mathcal{A}(Dv_q(t)+Dz_q(t))\|_{C^{0,1}}+l^{1/2-2\delta}\|\mathcal{A}(Dv_q+Dz_q)\|_{C_{[\frac{\sigma_q}{2},t]}^{1/2-2\delta} L^1}.\notag
\end{align}
Here $C^{0,1}$ is the Lipschitz-norm, i.e. $\|Q\|_{C^{0,1}}:=\sup_{x,y\in\mathbb{T}^d}\frac{|Q(x)-Q(y)|}{|x-y|}$ for all $Q:\mathbb{T}^d\to\mathbb{R}^{d\times d}$.

By \cite[Lemma 9]{DV} we have
\begin{align}
\|\nabla^2z^{in}(t)\|_{L^\infty}\lesssim(1+t^{-d/4-1})\|u_0\|_{L^2},\ \|\nabla z^{in}(t)\|_{L^\infty}\lesssim(1+t^{-d/4-1/2})\|u_0\|_{L^2}.\label{zin5}
\end{align}
Thus by Lemma \ref{lem:1}, (\ref{3.2}), (\ref{3.5})
we obtain
\begin{align}
l\|\mathcal{A}(Dv_q(t)+Dz_q(t))\|_{C^{0,1}}&\lesssim l\|\nabla^2v_q(t)+\nabla^2z_q(t)\|_{L^\infty}\lesssim l(\|v_q\|_{C_{t,x}^2}+\|\nabla^2z^{in}\|_{C_{[\frac{\sigma_q}{2},t]}L^\infty}+\|\nabla^2Z_q\|_{C_tL^\infty})\notag\\
&\lesssim l(\lambda_{q}^{\frac{3d}{2}+2}+\sigma_q^{-d/4-1}+\lambda_{q+1}^{\alpha/4})M_L^{1/2}\lesssim M_L^{1/2}\lambda_{q+1}^{-\alpha/8}\leq\frac{M_L\delta_{q+2}}{28},\notag
\end{align}
where we used the conditions $\alpha b>6d$ to deduce $\lambda_q^{\frac{3d}{2}}\leq\lambda_{q+1}^{\frac{\alpha}{4}}$, and used $l^{1/2}\leq\sigma_q^{d/4+1}$ to deduce $l\sigma_q^{-d/4-1}\leq l^{1/2}\leq \lambda_{q+1}^{-\alpha/8}$, and finally we used $\alpha>16\beta b$ and $a$ large enough to obtain $\lambda_{q+1}^{-\alpha/8+2\beta b}\leq \frac{1}{28}M_L^{1/2}.$

For the second term, by Lemma \ref{lem:1}, (\ref{3.2}), (\ref{3.5}) and (\ref{zin4})
we obtain
\begin{align}
l^{1/2-2\delta}&\|\mathcal{A}(Dv_q+Dz_q)\|_{C_{[\frac{\sigma_q}{2},t]}^{1/2-2\delta}L^1}\lesssim l^{1/2-2\delta}\|\nabla v_q+\nabla Z_q+\nabla z^{in}\|_{C_{[\frac{\sigma_q}{2},t]}^{1/2-2\delta} L^1}\notag\\
&\lesssim M_L^{1/2}l^{1/2-2\delta}(\lambda_{q}^{\frac{3d}{2}+2}+\lambda_{q+1}^{\frac\alpha4}+\sigma_{q}^{-1})\lesssim M_L^{1/2}\lambda_{q+1}^{-\alpha/8}\leq\frac{M_L\delta_{q+2}}{28},\notag
\end{align}
where we used the conditions $\alpha b>6d+8$ to deduce $\lambda_q^{\frac{3d}{2}+2}\leq\lambda_{q+1}^{\frac{\alpha}{4}}$, and used $l^{1/6}\leq\sigma_q$ to deduce $l^{1/2-2\delta}\sigma_q^{-1}\leq l^{1/3-2\delta}\leq \lambda_{q+1}^{-\alpha/8}$, and finally we used $\alpha>16\beta b$ and $a$ large enough to obtain $\lambda_{q+1}^{-\alpha/8+2\beta b}\leq \frac{1}{28}M_L^{1/2}.$

Similarly for $r\in(1,2],\nu_0=0$, by Lemma \ref{lem:1} we obtain
\begin{align}
\|{R}_{com21}(t)\|_{L^1}&\lesssim l^{r-1}\|\mathcal{A}(Dv_q(t)+Dz_q(t))\|_{C^{0,r-1}}+l^{(1/2-2\delta)(r-1)}\|\mathcal{A}(Dv_q+Dz_q)\|_{C_{[\frac{\sigma_q}{2},t]}^{(1/2-2\delta)(r-1)} L^1}.\notag\\
&\lesssim l^{r-1}\|\nabla^2 v_q(t)+\nabla^2 z_q(t)\|_{L^\infty}^{r-1}+l^{(1/2-2\delta)(r-1)}\|\nabla v_q+\nabla z_q\|_{C_{[\frac{\sigma_q}{2},t]}^{1/2-2\delta} L^1}^{r-1}.\notag\\
&\lesssim M_L^{1/2}\lambda_{q+1}^{-\alpha(r-1)/8}\leq\frac{M_L\delta_{q+2}}{14},\notag
\end{align}
where $C^{0,r-1}$ equals to the H\"older-norm $C^{r-1}$ for $r<2$, and equals to the Lipschitz-norm $C^{0,1}$ for $r=2$. Then we used the conditions on the parameters to deduce $\alpha(r-1)>16\beta b$ and $a$ large enough such that $\lambda_{q+1}^{-\alpha(r-1)/8+2\beta b}\leq\frac{1}{14}M_L^{1/2}$.

The case $ r\in(2,\frac{3d+2}{d+2})$ follows from
\begin{align*}
\|{R}_{com21}(t)\|_{L^1}&\lesssim l\|\mathcal{A}(Dv_q(t)+Dz_q(t))\|_{C^{0,1}}+ l^{1/2-2\delta}\|\mathcal{A}(Dv_q+Dz_q)\|_{C_{[\frac{\sigma_q}{2},t]}^{1/2-2\delta} L^1}\\
&\lesssim (l\|\nabla^2v_q(t)+\nabla^2z_q(t)\|_{L^\infty}+l^{1/2-2\delta}\|\nabla v_q+\nabla z_q\|_{C_{[\frac{\sigma_q}{2},t]}^{1/2-2\delta} L^1})(1+\|\nabla v_q+\nabla z_q\|^{r-2}_{C_{[\frac{\sigma_q}{2},t]}L^\infty}).\notag
\end{align*}
The estimate of the first part is the same as the case $r\in(1,2],\nu_0>0$ and the estimate of the second part follows from (\ref{3.2}), (\ref{3.4}) and (\ref{zin5}). Thus we obtain
\begin{align*}
\|{R}_{com21}(t)\|_{L^1}&\lesssim(l\sigma_q^{-d/4-1}+l^{1/2-2\delta}\lambda_{q}^{\frac{3d}{2}+2}+l^{1/2-2\delta}\lambda_{q+1}^{\frac\alpha4}+l^{1/2-2\delta}\sigma_{q}^{-1})(\lambda_q^{d+1}+\lambda_{q+1}^{\frac\alpha4}+\sigma_{q}^{-d/4-1/2})M_L\notag\\
&\lesssim l^{1/2-2\delta}(\lambda_{q+1}^{\frac\alpha4}+\sigma_{q}^{-1})(\lambda_{q+1}^{\frac\alpha4}+\sigma_{q}^{-d/4-1/2})M_L\notag\\
&\lesssim M_L\lambda_{q+1}^{-\alpha/8}\leq\frac{M_L\delta_{q+2}}{14},\notag
\end{align*}
where in the second inequality we used the conditions $\alpha b>6d+8$ to deduce $\lambda_q^{\frac{3d}{2}+2}\leq\lambda_{q+1}^{\frac{\alpha}{4}}$ and used $l^{1/2}\leq\sigma_q^{d/4}$ to deduce $l\sigma_q^{-d/4-1}\leq l^{1/2-2\delta}\sigma_q^{-1}$, and in the third inequality we used the condition $l^{1/6}\leq\sigma_q^{d/4+1/2}$ to deduce $l^{1/2-2\delta}\lambda_{q+1}^{\frac\alpha4}\sigma_{q}^{-d/4-1/2}\leq l^{1/4}\lambda_{q+1}^{\frac\alpha4}\leq\lambda_{q+1}^{-\alpha/8}$, $l^{1/2-2\delta}\sigma_q^{-1-d/4-1/2}\leq \lambda_{q+1}^{-\alpha/8}$ and $l^{1/2-2\delta}\lambda_{q+1}^{\alpha/2}\leq\lambda_{q+1}^{-\alpha/8}$. Finally we used $\alpha>16\beta b$ and $a$ large enough to obtain $\lambda_{q+1}^{-\alpha/8+2\beta b}\leq \frac{1}{14}.$

Similarly to the estimate of nonlinear error, by (\ref{3.2}) and (\ref{bdd:Dz1}),  we obtain for $ r\in(1,2],\nu_0>0$
\begin{align}
\|{R}_{com22}(t)\|_{L^1}&\lesssim\|\nabla v_l(t)-\nabla v_q(t)+\nabla z_l(t)-\nabla z_q(t)\|_{L^1}\lesssim l\|\nabla v_q\|_{C_{x,t}^1}+\|\nabla z_l(t)-\nabla z_q(t)\|_{L^1}.\notag\\
&\lesssim M_L^{1/2}l\lambda_{q}^{\frac{3d}{2}+2}+M_L^{1/2}\lambda_{q+1}^{-\alpha/8}\lesssim M_L^{1/2}\lambda_{q+1}^{-\alpha/8}\leq\frac{M_L\delta_{q+2}}{14},\notag
\end{align}
where we used the conditions $\alpha b>6d$ to deduce $\lambda_q^{\frac{3d}{2}}\leq\lambda_{q+1}^{\frac{\alpha}{4}}$ and used $\alpha>16\beta b$ and $a$ large enough to obtain $\lambda_{q+1}^{-\alpha/8+2\beta b}\leq \frac{1}{14}M_L^{1/2}.$

For $ r\in(1,2],\nu_0=0$ similarly we obtain
\begin{align}
\|{R}_{com22}(t)\|_{L^1}\lesssim\|\nabla v_l(t)-\nabla v_q(t)+\nabla z_l(t)-\nabla z_q(t)\|^{r-1}_{L^1}
\lesssim M_L^{1/2}\lambda_{q+1}^{-\alpha(r-1)/8}\leq\frac{M_L\delta_{q+2}}{14},\notag
\end{align}
where we used the conditions on the parameters to deduce $\alpha(r-1)>16\beta b$ and $a$ large enough such that $\lambda_{q+1}^{-\alpha(r-1)/8+2\beta b}\leq\frac{1}{14}M_L^{1/2}$.

For $r\in(2,\frac{3d+2}{d+2})$, by (\ref{3.2}), (\ref{3.4}), (\ref{3.5}) and (\ref{bdd:Dz1})
\begin{align}
\|{R}_{com22}\|_{C_tL^1}&\lesssim(\|\nabla v_l-\nabla v_q\|_{C_{t}L^{r-1}}+\|\nabla z_l(t)-\nabla z_q(t)\|_{L^{r-1}})\notag\\
&\ \ \ \ \ \ \ \ \times(1+\|\nabla v_l+\nabla z_l\|_{C_{t}L^{r-1}}^{r-2}+\|\nabla v_q+\nabla z_q\|
_{C_{t}L^{r-1}}^{r-2})\notag\\
&\lesssim (l\|v_q\|_{C_{t,x}^2}+\|\nabla z_l(t)-\nabla z_q(t)\|_{L^{r-1}})(1+\|\nabla z_q\|
_{C_{t}L^{r-1}}^{r-2}+\|\nabla v_q\|
_{C_{t}L^{r-1}}^{r-2})\notag\\
&\lesssim M_L^{1/2}(l\lambda_{q}^{\frac{3d}{2}+2}+l^{1/4}+l^{1/2-2\delta}\lambda_{q+1}^{\frac\alpha4})(1+\lambda_{q}^{d+1}+\lambda_{q+1}^{\alpha/4})M_L^{1/2}\notag\\
&\lesssim M_L(l^{1/4}\lambda_{q+1}^{\frac{\alpha}4}+l^{1/2-2\delta}\lambda_{q+1}^{\frac\alpha2})\lesssim M_L\lambda_{q+1}^{-\alpha/8}\leq\frac{M_L\delta_{q+2}}{14},\notag
\end{align}
where we used the conditions $\alpha b>6d+8$ to deduce $\lambda_q^{\frac{3d}{2}+2}\leq\lambda_{q+1}^{\frac{\alpha}{4}}$ and used $\alpha>16\beta b$ and $a$ large enough to obtain $\lambda_{q+1}^{-\alpha/8+2\beta b}\leq \frac{1}{14}.$

We continue with $R_{com3}$. By (\ref{3.2}), (\ref{3.3}), (\ref{zin11}) and (\ref{zin13}) we have
\begin{align}
\|R_{com3}(t)\|_{L^1}&\lesssim M_0(M_L+K+qA)^{1/2}(\|z_l(t)-z_{q}(t)\|_{L^2}+\|z_q(t)-z_{q+1}(t)\|_{L^2})\notag\\
&\lesssim M_0(M_L+K+qA)^{1/2}M_L^{1/2}(l^{1/2-2\delta}+l^{1/2}\sigma_{q}^{-1/2}+\lambda_{q+1}^{-\frac{\alpha(1+\sigma)}{2d+4}})\notag\\
&\lesssim M_0(M_L+K+qA)(\lambda_{q+1}^{-\alpha/8}+\lambda_{q+1}^{-\frac{\alpha(1+\sigma)}{2d+4}})\leq\frac{M_L\delta_{q+2}}{7},\label{3.21}
\end{align}
where we used (\ref{a2}) and the conditions $l^{1/4}\leq \sigma_q^{1/2}$ to deduce $l^{1/2}\sigma_{q}^{-1/2}\leq l^{1/4}$.

Summarizing estimates (\ref{3.15})-(\ref{3.21}) we obtain
$$\|\mathring{R}_{q+1}(t)\|_{L^1}\leq M_L\delta_{q+2}.$$

\textbf{2} Let $t\in(\frac{\sigma_q}2\wedge T_L,\sigma_q\wedge T_L]$. If $T_L\leq \frac{\sigma_q}2$
then there is nothing to estimate, hence we may assume $\frac{\sigma_q}2 < T_L$ and $t\in(\frac{\sigma_q}2, \sigma_q \wedge T_L]$.
Then we decompose $\mathring{R}_l =\chi^2\mathring{R}_l + (1-\chi^2)\mathring{R}_l $ The first
part $\chi^2\mathring{R}_l $ is canceled (up to the oscillation error $\chi^2R_{osc}$) by $\tilde{w}^{(p)}_{ q+1}\otimes \tilde{w}^{(p)}_{ q+1}= \chi^2w^{(p)}_{ q+1}\otimes w^{(p)}_{ q+1}$ and
$\chi^2\partial_tw^{(t)}_{ q+1} = \partial_t\tilde{w}^{(t)}_{q+1}-(\chi^2)' w^{(t)}_{ q+1}$. So in this case in the definition of $\mathring{R}_{q+1}$ most terms are similar to in the case \textbf{1}. and can be estimated similarly.
We only have to consider $(1-\chi^2)\mathring{R}_l$ and
\begin{align}
\div\mathring{R}_{cut}:=\chi'(t)(w_{q+1}^{(p)}(t)+w_{q+1}^{(c)}(t))+(\chi^2)'(t)w_{q+1}^{(t)}(t).\notag
\end{align}
We know
$$\|(1-\chi^2)\mathring{R}_l(t)\|_{L^1}\leq \sup_{s\in[t-\sigma_q,t]}\|\mathring{R}_q(s)\|_{L^1}.$$
As for $\mathring{R}_{cut}$ by (\ref{2.18+-})-(\ref{2.18**}) we have
\begin{align*}
\|\mathring{R}_{cut}(t)\|_{L^1}&\leq\|\chi'(t)(w_{q+1}^{(p)}+w_{q+1}^{(c)})\|_{L^{1+\epsilon}}+\|(\chi^2)'w_{q+1}^{(t)}\|_{L^{1+\epsilon}}\\
&\lesssim \sigma_q^{-1}(\|w_{q+1}^{(p)}\|_{L^{1+\epsilon}}+\|w_{q+1}^{(c)}\|_{L^{1+\epsilon}}+\|w_{q+1}^{(t)}\|_{L^{1+\epsilon}})\\
&\leq \sigma_q^{-1}(M_L+qA+\gamma_{q+1})\lambda_{q+1}^{(4d+4)\alpha+(\frac{1}{1+\epsilon}-\frac{1}{2})(\frac{2d+1}{N}-1)}\\
&\lesssim\lambda_{q+1}^{(4d+8)\alpha-\frac{d}{4}+{d\epsilon}}\leq\frac{M_L\delta_{q+2}}{14},
\end{align*}
where we choose $\epsilon>0$ small enough such that ${d\epsilon}<\alpha$ and we used the conditions $\sigma_q^{-1}\leq l^{-1}$ and $\frac{d}{4}-(4d+9)\alpha>2\beta b$.

For $R_{com}$, $R_{com1}$ and $R_{com2}$ since $t >\sigma_q/2$ and $l\leq\frac{\sigma_q}4$ we have a similar bound
as in the first case.

\textbf{3} For $t\in[0,\frac{\sigma_q}2\wedge T_L]$ we know $v_l(t)=v_{q+1}(t)=0$ and
\begin{align}
\mathring{R}_{q+1}=&\mathring{R}_l+(\mathcal{A}(Dz_{q})*_x\phi_l)*\varphi_l-\mathcal{A}(Dz_{q+1})\notag\\
&+2D z_{q+1}-2D z_l+z_{q+1}\mathring{\otimes}z_{q+1}-((z_q\mathring{\otimes}z_q)*_x\phi_l)*_t\varphi_l.\notag
\end{align}
Similarly as (\ref{adz0}), by Lemma \ref{lem:1} and the choice of $M_L$ we have $\|{\mathcal{A}(Dz_{q+1})}\|_{L^1}\leq \frac{1}{2}M_L,\|{\mathcal{A}(Dz_{q})}\|_{L^1}\leq \frac{1}{2}M_L$. Thus
\begin{align}
\|\mathring{R}_{q+1}(t)\|_{L^1}&\leq\sup_{s\in[t-l,t]}\|\mathring{R}_{q}(s)\|_{L^1}+M_L+2\|\nabla z_{q+1}\|_{L^1}+2\|\nabla z_{q}\|_{L^1}+2\|z\|_{L^2}^2\notag\\
&\leq\sup_{s\in[t-l,t]}\|\mathring{R}_{q}(s)\|_{L^1}+M_L+4(2\pi)^{\frac d2}(\|\nabla Z\|_{L^2}+\|\nabla u_0\|_{L^{r^*}})+2(\|Z\|_{L^2}+\|u_0\|_{L^2})^2\notag\\
&\leq\sup_{s\in[t-l,t]}\|\mathring{R}_{q}(s)\|_{L^1}+M_L+4(2\pi)^{\frac d2}(L+M)+2(L+M)^2\notag\\
&\leq\sup_{s\in[t-l,t]}\|\mathring{R}_{q}(s)\|_{L^1}+2M_L,\notag
\end{align}
which implies (\ref{3.9}) and (\ref{3.6}), (\ref{3.7}) at the level $q+1$.

We finish the proof of Proposition \ref{prop:1}.

\section{Uniqueness in 3D}\label{tctr}
In this section we prove Theorem \ref{thm:d=3} by extending the classical results \cite[Chapter 5 Theorem 4.37]{MNRR} in the PDE case to the stochastic setting. By \cite[Theorem 2.2.1]{YN}  we only need to consider the case $\frac{11}{5}\leq r<\frac{5}{2}$. In the following we consider the equation on $[0,T]$ for simplicity.
\bt\label{prop:rr}
Let $\frac{11}{5}\leq r<\frac{5}{2}$, $\nu_0>0$, and $G$ satisfying $\|\nabla G\|_{L_2(U,L^2_\sigma)}<\infty$, $T>0$. Then there exists a unique probabilistically strong solution $u$ in $L^r([0,T];W^{1,r})\cap C([0,T];L_\sigma^2)$ to (\ref{eq:pl}) with the initial value $u_0\in H^1$. Moreover, the solution satisfies $\sup_{t\in[0,T]}\|\nabla u(t)\|_{L^2}<\infty$ $\mathbf{P}$-a.s.
\et
\begin{proof}
Our proof is divided into the following four steps:

(a). We consider the following Galerkin approximation
\begin{align}\dif u^N=\mathbb{P}_N\mathbb{P}\div\mathcal{A}(Du^N)\dif t&-\mathbb{P}_N\mathbb{P}(u^N\cdot \nabla)u^N\dif t+\mathbb{P}_N\dif B,\notag\\\div u^N&=0,\\u^N(0)&=\mathbb{P}_Nu_0.\notag\end{align}Here $\mathbb{P}_N$ is the finite dimensional projection and $\mathbb{P}$ is the Leray Projection.. By the proof of \cite[Theorem 2.1.3]{YN} the law of $u^N$ is tight in $\mathbf{S}:= C([0,T];W^{-\zeta,\frac{r}{r-1}}) \cap L^{\tilde{r}}([0,T]; W^{1,\tilde{r}})$, for some $\zeta>0, \tilde{r}\in(1,r)$. Then by the Skorohod's theorem there exists a probability space $(\tilde{\Omega}, \tilde{\mathcal{F}}, \tilde{\bP})$ and $\mathbf{S}$-valued random variables $\tilde{u}^N$ and ${u}$ such that

(1).\ $\tilde{u}^N$ has the same law as $u^N$ for each $n\in\mathbb{N}$;

(2).\ $\tilde{u}^N\to u$ in $\mathbf{S}$ $\tilde{\bP}$-a.s.

Then by the proof of \cite[Theorem 2.1.3]{YN} we obtain $u$ is a probabilistically weak solution to (\ref{eq:pl}) satisfying
 \begin{align}
\mathbb{E}[\sup_{t\in[0,T]}\|u(t)\|_{L^2}^2+\int_0^T\|\nabla u(t)\|_{L^r}^r\dif t]<\infty.\label{etnablavrr}
\end{align}

(b). The solution $u$ satisfies $\sup_{t\in[0,T]}\|\nabla u(t)\|_{L^2}<\infty$ $\tilde{\bP}$-a.s.

We denote $\lambda=\frac{6-2r}{3r-5}$ which is smaller than $1$ by $r\geq\frac{11}{5}$. By the proof in \cite[Lemma 3.2.3]{YN}
we obtain for $\lambda<1$ the Galerkin approximation $u^N$ satisfies
\begin{align}
\frac{1}{1-\lambda}&(1+\|\nabla u^N(t)\|_{L^2}^2)^{1-\lambda}-\frac{1}{1-\lambda}(1+\|\nabla \mathbb{P}_Nu_0\|_{L^2}^2)^{1-\lambda} \notag\\
&\lesssim \int_0^t (1+\|\nabla u^N\|_{L^r})^r\dif s+\|\nabla G\|_{L_2}^2t+\int_0^t \frac{\langle -\Delta u^N,\dif \mathbb{P}_NB\rangle}{(1+\|\nabla u^N(t)\|_{L^2}^2)^\lambda}.\label{123456}
\end{align}
By Burkholder-Davis-Gundy's inequality and Young's inequality we obtain
\begin{align*}
\mathbb{E}[\sup_{t\in[0,T]}|\int_0^t \frac{\langle -\Delta u^N,\dif \mathbb{P}_NB\rangle}{(1+\|\nabla u^N(t)\|_{L^2}^2)^\lambda}]
&\lesssim\mathbb{E}[\int_0^T\frac{\|\nabla u^N(t)\|_{L^2}^2\|\nabla G\|_{L_2}^2}{(1+\|\nabla u^N(t)\|_{L^2}^2)^{2\lambda}}\dif t]^{1/2}\\
&\leq \epsilon\mathbb{E}[\sup_{t\in[0,T]}(1+\|\nabla u^N(t)\|_{L^2}^2)^{1-\lambda}]+C\|\nabla G\|_{L_2}^2.
\end{align*}
Then
taking supremum over $t$ and taking expectations in (\ref{123456}) and by (\ref{etnablavrr}) we obtain
\begin{align*}
\mathbb{E}[\sup_{t\in[0,T]}(1+\|\nabla u^N(t)\|_{L^2}^2)^{1-\lambda}]<\infty.
\end{align*}
Here the constant is independent of $N$.

For the case $\lambda=1$ similarly we obtain
\begin{align*}
\mathbb{E}[\sup_{t\in[0,T]}\ln(1+\|\nabla u^N(t)\|_{L^2}^2)]<\infty.
\end{align*}

Thus we obtain $\sup_{t\in[0,T]}\|\nabla \tilde{u}^N(t)\|_{L^2}<\infty$ a.s. and the constant is independent of $N$. Letting $N\to\infty$ we conclude the proof by lower semi-continuity.

(c). $u\in C([0,T];L_\sigma^2)$ $\tilde{\bP}$-a.s.

To prove the continuity in $L_\sigma^2$, by (\ref{etnablavrr}) and continuity in $W^{-\zeta,\frac{r}{r-1}}$ the process $u$ is weakly continuous in $L^2_\sigma$. Therefore we only need to prove $t\to\|u(t)\|_{L^2}$ is continuous.

By It\^o's formula we obtain
\begin{align}
\|u(t)\|_{L^2}^2&=-2\int_0^t\int\mathcal{A}(Du):Du\dif x\dif s+\|G\|_{L_2}^2t+2\int_0^t\langle u,\dif B\rangle.\notag
\end{align}
Here we used  \cite[Chapter 5 Lemma 2.44]{MNRR} to deduce the nonlinear term $\int_0^t\langle (u\cdot\nabla)u,u \rangle\dif s$ is finite and equals zero for $r\geq \frac{11}{5}$ and $u\in L^r([0,T];W^{1,r})\cap L^\infty([0,T];L_\sigma^2)$.

Notice that the third term on the right side is a continuous martingale so we only need to consider the first term on the right side. By the definition of $\mathcal{A}$ we obtain for $0\leq t_1<t_2\leq T$
$$|\int_{t_1}^{t_2}\int\mathcal{A}(Du):Du\dif x\dif s|\lesssim \int_{t_1}^{t_2} (\|\nabla u\|_{L^{r}}^r+1)\dif s.$$
Together with (\ref{etnablavrr}) we get the continuity of the first term by the dominated convergence Theorem.

(d). There exists a unique probabilistically strong solution $u$ to (\ref{eq:pl}) in $L^r([0,T];W^{1,r})\cap C([0,T];L_\sigma^2)$.

First
 assume that there exists two weak solutions $u,v\in L^r([0,T];W^{1,r})\cap C([0,T];L_\sigma^2)$ on the same stochastic basis $(\Omega_1, \mathcal{F}_1,\mathbf{P}_1)$ satisfying $\sup_{t\in[0,T]}\|\nabla u(t)\|_{L^2}+\sup_{t\in[0,T]}\|\nabla v(t)\|_{L^2}<\infty$ $\mathbf{P}_1$-a.s.  Set $w=u-v$ then we obtain
\begin{align*}
\dif w=\mathbb{P}[\div{\mathcal{A}(Du)}&-\div{\mathcal{A}(Dv)}]\dif t-\mathbb{P}[(u\cdot \nabla)u-(v\cdot \nabla)v]\dif t,\\
\div w&=0,\\
w(0)&=0.
\end{align*}
Here $\mathbb{P}$ is the Leray Projection.

Taking inner product in $L^2$ with $w$ we get that for all $t\in[0,T]$
\begin{align}
\|w(t)\|_{L^2}^2&=-2\int_0^t\int(\mathcal{A}(Du)-\mathcal{A}(Dv))(Dw)\dif x\dif s-2\int_0^t\langle w\cdot \nabla u,w\rangle \dif s.\notag
\end{align}
Here we used \cite[Chapter 5 Lemma 2.44]{MNRR} to deduce the nonlinear term $\int_0^t\langle (u\cdot\nabla)w,w \rangle\dif s$ is finite and equals zero for $r\geq \frac{11}{5}$ and $u,v\in L^r([0,T];W^{1,r})\cap L^\infty([0,T];L_\sigma^2)$.
By \cite[Chapter 5 Lemma 1.19]{MNRR} and Korn’s inequality (\cite[Section 5.1.1]{MNRR}) we obtain
$$\int_0^t\|\nabla w(s)\|_{L^2}^2\dif s\lesssim\int_0^t\int(\mathcal{A}(Du)-\mathcal{A}(Dv))(Dw)\dif x\dif s,$$
which together with H\"{o}lder's inequality implies
$$\|w(t)\|_{L^2}^2+C\int_0^t\|\nabla w(s)\|_{L^2}^2\dif s\lesssim\int_0^t\|\nabla u(s)\|_{L^2}\|w(s)\|_{L^4}^2\dif s.$$
Moreover, using Sobolev's embedding and the interpolation inequality we get
$$\|w(s)\|_{L^4}\lesssim \|w(s)\|_{L^2}^{1/4}\|\nabla w(s)\|_{L^2}^{3/4},$$
which by Young's inequality implies that
\begin{align*}
\|w(t)\|_{L^2}^2+C\int_0^t\|\nabla w(s)\|_{L^2}^2\dif s&\lesssim\int_0^t\|\nabla u(s)\|_{L^2}\|w(s)\|_{L^2}^{1/2}\|\nabla w(s)\|_{L^2}^{3/2}\dif s\\
&\leq C\int_0^t\|\nabla w(s)\|_{L^2}^2\dif s+C'\int_0^t\|\nabla u(s)\|_{L^2}^4\|w(s)\|_{L^2}^2\dif s,
\end{align*}
i.e.
\begin{align*}
\|w(t)\|_{L^2}^2\leq C'\int_0^t\|\nabla u(s)\|_{L^2}^4\|w(s)\|_{L^2}^2\dif s.
\end{align*}
Hence using the condition $\sup_{t\in[0,T]}\|\nabla u(t)\|_{L^2}<\infty$ and Gronwall's lemma we  obtain the pathwise uniqueness for the solutions satisfying $\sup_{t\in[0,T]}\|\nabla u(t)\|_{L^2}<\infty$. Thus by the Yamada-Watanabe Theorem we obtain there exists a unique probabilistically strong solution $u$ satisfying
$\sup_{t\in[0,T]}\|\nabla u(t)\|_{L^2}<\infty$.

Now we prove the solution is unique in $L^r([0,T];W^{1,r})\cap C([0,T];L_\sigma^2)$. In fact, for a probabilistically strong solution $u^*\in L^r([0,T];W^{1,r})\cap C([0,T];L_\sigma^2)$, by the same argument as above we obtain $u=u^*$, which implies uniqueness.
\end{proof}

\appendix
 \renewcommand{\appendixname}{Appendix~\Alph{section}}
  \renewcommand{\theequation}{A.\arabic{equation}}

  \section{Estimates of $\rho$ and $a_{(\xi)}$}
  \label{s:appA}

We first give the estimates of $\rho$ and $a_{(\xi)}$ for the case $t\in[0,T_L]$. By the definition of $\rho$ we have
$$\|\rho\|_{C_tL^p}\leq 2l(2\pi)^{d/p}+2\|\mathring{R}_l\|_{C_tL^p}+\gamma_{q+1}(2\pi)^{d/p-d}.$$
Furthermore, by mollification estimates, the embedding $W^{d+1,1}\subset L^\infty$ and (\ref{3.7}) we obtain for $N\geq 0,t\in[0,T_L]$
$$\|\mathring{R}_l\|_{C_{t,x}^N}\lesssim l^{-d-1-N}(M_L+qA),$$
which in particular leads to
$$\|\rho\|_{C_{t,x}^0}\lesssim l+l^{-d-1}(M_L+qA)+\gamma_{q+1}\lesssim l^{-d-1}(M_L+qA)+\gamma_{q+1}.$$
Next we estimate the $C_{t,x}^N$-norm for $N\in\mathbb{N}$. We apply the chain rule in \cite[Proposition C.1]{TCP}  to $f(z)=\sqrt{l^2+z^2},|D^mf(z)|\lesssim l^{-m+1}$ to obtain
\begin{align*}
\|\sqrt{l^2+|\mathring{R}_l|^2}\|_{C_{t,x}^N}&\lesssim\|\sqrt{l^2+|\mathring{R}_l|^2}\|_{C_{t,x}^0}+\|Df\|_{C^0}\|\mathring{R}_l\|_{C_{t,x}^N}+\|Df\|_{C^{N-1}}\|\mathring{R}_l\|_{C_{t,x}^1}^N\\
&\lesssim l^{-d-1-N}(M_L+qA)+l^{-N+1}l^{-(d+2)N}(M_L+qA)^N,
\end{align*}
then by $M_L+qA\lesssim l^{-1}$ we deduce for $N\geq1$
\begin{align}\label{A.1}
\|\rho\|_{C_{t,x}^N}&\lesssim\|\sqrt{l^2+|\mathring{R}_l|^2}\|_{C_{t,x}^N}+\gamma_{q+1}
\notag\\
&\lesssim(l^{-d-1-N}+l^{-N+1}l^{-(d+2)N}l^{-N+1})(M_L+qA)+\gamma_{q+1}\notag\\
&\lesssim l^{2-(d+4)N}(M_L+qA)+\gamma_{q+1}.
\end{align}
Next we estimate $a_{(\xi)}$. By (\ref{3.13*}) and Lemma \ref{lem:3},
\begin{align}
\|a_{(\xi)}\|_{C_tL^2}&\leq \|\rho\|_{C_tL^1}^{1/2}\|\gamma_\xi\|_{C^0(B_{1/2}(\Id))}\notag\\
&\leq\frac{M^*}{8|\Lambda|(1+(2\pi)^d)^{1/2}}(2l(2\pi)^d+2M_L+2qA+\gamma_{q+1})^{1/2}\notag\\
&\leq \frac{M^*}{4|\Lambda|}((M_L+qA)^{1/2}+\gamma_{q+1}^{1/2}).
\end{align}
Let us now estimate the $C_{t,x}^N$-norm. By Leibniz rule we get
$$\|a_{(\xi)}\|_{C_{t,x}^N}\lesssim\sum_{m=0}^N\|\rho^{1/2}\|_{C_{t,x}^m}\|\gamma_\xi(\Id-\frac{\mathring{R}_l}{\rho})\|_{C_{t,x}^{N-m}}.$$
Apply \cite[Proposition C.1]{TCP} to $f(z)=z^{1/2}$, $|D^mf(z)|\lesssim|z|^{1/2-m}$, for $m=1,...,N$, and using (\ref{A.1}) we obtain
\begin{align*}
\|\rho^{1/2}\|_{C_{t,x}^m}&\lesssim \|\rho^{1/2}\|_{C_{t,x}^0}+l^{-1/2}\|\rho\|_{C_{t,x}^m}+l^{1/2-m}\|\rho\|_{C_{t,x}^1}^m\\
    &\lesssim (l^{-(d+1)/2}+l^{1-(d+4)m})((M_L+qA)^{1/2}+\gamma_{q+1}^{1/2})+l^{1/2-m}l^{-(d+2)m}(M_L+qA+\gamma_{q+1})^{m}\\
&\lesssim l^{1-(d+4)m}((M_L+qA)^{1/2}+\gamma_{q+1}^{1/2}),
\end{align*}
where we used $M_L+qA+K\leq l^{-1}$.

Next we estimate $\gamma_\xi(\Id-\frac{\mathring{R}_l}{\rho})$, by \cite[Proposition C.1]{TCP} we need to estimate
$$\|\frac{\mathring{R}_l}{\rho}\|_{C_{t,x}^{N-m}}+\|\frac{\nabla_{t,x}\mathring{R}_l}{\rho}\|_{C_{t,x}^0}^{N-m}+\|\frac{\mathring{R}_l}{\rho^2}\|_{C_{t,x}^0}^{N-m}\|\rho\|_{C_{t,x}^1}^{N-m}.$$
We use $\rho\geq l$ and $M_L+qA\leq l^{-1}$ to have
$$\|\frac{\nabla_{t,x}\mathring{R}_l}{\rho}\|_{C_{t,x}^0}^{N-m}\lesssim l^{-N+m}l^{-(d+2)(N-m)}(M_L+qA)^{N-m}\lesssim l^{-(d+4)(N-m)},$$
and in view of $|\frac{\mathring{R}_l}{\rho}|\leq 1$
$$\|\frac{\mathring{R}_l}{\rho^2}\|_{C_{t,x}^0}^{N-m}\lesssim\|\frac{1}{\rho}\|_{C_{t,x}^0}^{N-m}\lesssim l^{-N+m},$$
and by (\ref{A.1}) and $M_L+qA+K\leq l^{-1}$
$$\|\rho\|_{C_{t,x}^1}^{N-m}\lesssim l^{-(d+2)(N-m)}(M_L+qA+\gamma_{q+1})^{-(N-m)}\lesssim l^{-(d+3)(N-m)}.$$
Moreover, we write
\begin{align*}
\|\frac{\mathring{R}_l}{\rho}\|_{C_{t,x}^{N-m}}&\lesssim  \sum_{k=0}^{N-m}\|\mathring{R}_l\|_{C_{t,x}^k}\|\frac{1}{\rho}\|_{C_{t,x}^{N-m-k}}.
\end{align*}
Using (\ref{A.1}) and $M_L+qA+K\leq l^{-1}$
\begin{align*}
\|\frac{1}{\rho}\|_{C_{t,x}^{N-m-k}}\lesssim \|\frac{1}{\rho}\|_{C_{t,x}^{0}}+l^{-2}\|{\rho}\|_{C_{t,x}^{N-m-k}}+l^{-N+m+k-1}\|{\rho}\|_{C_{t,x}^1}^{N-m-k}\\
\lesssim l^{-2}l^{2-(d+4)(N-m-k)-1}+l^{-(N-m-k)-1}l^{-(d+2)(N-m-k)-(N-m-k)}\lesssim l^{-(d+4)(N-m-k)-1}.
\end{align*}
Thus we obtain
\begin{align*}
\|\frac{\mathring{R}_l}{\rho}\|_{C_{t,x}^{N-m}}&\lesssim  \sum_{k=0}^{N-m-1}l^{-d-2-k}l^{-(d+4)(N-m-k)-1}+l^{-d-2-(N-m)}l^{-1}\lesssim l^{-(d+3)-(d+4)(N-m)}.
\end{align*}
Finally the above bounds lead to
$$\|\gamma_\xi(\Id-\frac{\mathring{R}_l}{\rho})\|_{C_{t,x}^{N-m}}\lesssim l^{-(d+3)-(d+4)(N-m)}+l^{-(d+4)(N-m)}\lesssim l^{-(d+3)-(d+4)(N-m)}.$$
Combining this with the bounds for $\rho^{1/2}$ above yields for $N\in\mathbb{N}$
\begin{align*}
\|a_{(\xi)}\|_{C_{t,x}^N}&\lesssim(l^{-(d+1)/2}l^{-(d+3)-(d+4)N}+\sum_{m=1}^{N-1}l^{1-(d+4)m}l^{-(d+3)-(d+4)(N-m)}+l^{1-(d+4)N})((M_L+qA)^{1/2}+\gamma_{q+1}^{1/2})\\
&\lesssim l^{-2d-2-(d+4)N}((M_L+qA)^{1/2}+\gamma_{q+1}^{1/2}),
\end{align*}
where the final bound is also valid for $N=0$.

For $t\in(4\sigma_q\wedge T_L,T_L]$, using $l\leq \sigma_q$ and (\ref{3.6}), the estimate is similar, then we have (\ref{al21}) and (\ref{3.12}).

  \renewcommand{\theequation}{B.\arabic{equation}}
  \section{Some technical tools}\label{tec}
  \subsection{Control of $\mathcal{A}$}\label{coa}
We collect the growth estimate for $\mathcal{A}(Q)$ required in the following.
\bl\label{lem:1}$($\cite[Lemma 1]{BMS}$)$
Let $\mathcal{A}(Q)=(\nu_0+\nu_1|Q|)^{r-2}Q$ with $
\nu_0\geq0,\nu_1\geq 0, r\in(1,\infty)$. Then
\begin{align}
|\mathcal{A}(Q)-\mathcal{A}(P)|\leq\begin{cases}
C_{\nu_1}|Q-P|^{r-1},\ \ &for\ \nu_0=0,r\leq2,\\
C_{\nu_0}|Q-P|,\ \ \ \ \ \ &for\ \nu_0>0,r\leq2,\\
C_{r,\nu_0,\nu_1}|Q-P|(1+|P|^{r-2}+|Q|^{r-2}),\ \ &for\ r> 2.
\end{cases}\label{A1}
\end{align}

\el
\subsection{Antidivergence operators}\label{tamr}
We first recall the following antidivergence operator $\mathcal{R}$ as in \cite[Appendix B.2]{AX}, which acts on vector fields $v$ with $\int_{\mathbb{T}^d} v\dif x = 0$ as
$$(\mathcal{R}v)_{ij}=\mathcal{R}_{ijk}v_k,$$
where

$$\mathcal{R}_{i j k}=\frac{2-d}{d-1} \Delta^{-2} \partial_{i} \partial_{j} \partial_{k}-\frac{1}{d-1} \Delta^{-1} \partial_{k} \delta_{i j}+\Delta^{-1} \partial_{i} \delta_{j k}+\Delta^{-1} \partial_{j} \delta_{i k}.$$
Then $\mathcal{R}v(x)$ is a symmetric trace-free matrix for each $x \in \mathbb{T}^d$, and $\mathcal{R}$ is a right
inverse of the $\div$ operator, i.e. $\div(\mathcal{R}v) = v$.
By a direct computation we have for any divergence-free $v\in C^\infty(\mathbb{T}^d,\mathbb{R}^d)$
\begin{align}\label{rdeltav}
\mathcal{R}\Delta v=\nabla v+\nabla^{T} v=2Dv.
\end{align}
We can show that $\mathcal{R}$ is bounded on $L^p(\mathbb{T}^d)$ for any $1\leq p \leq\infty$. In the following we use $C^\infty(\mathbb{T}^d,\mathbb{R}^d)$ to denote the space of smooth functions from $\mathbb{T}^d$ to $\mathbb{R}^d$. We use $C_0^\infty(\mathbb{T}^d,\mathbb{R}^d)$ to denote the subspace of functions with zero spatial mean. Similarly we define $C^\infty(\mathbb{T}^d,\mathbb{R}^{d\times d})$ and $C_0^\infty(\mathbb{T}^d,\mathbb{R}^{d\times d})$.
\bt$($\cite[Theorem B.3]{AX}$)$\label{bb}
 Let $1\leq p \leq\infty$. For any vector field $f \in C^\infty(\mathbb{T}^d, \mathbb{R}^d)$, there holds
$$\|\mathcal{R}f\|_{ L^p}\lesssim \| f\|_{ L^p}.$$
In particular, if $f \in C_0^\infty(\mathbb{T}^d, \mathbb{R}^d)$, then
$$\|\mathcal{R}f(\sigma\cdot)\|_{ L^p}\lesssim\sigma^{-1} \| f\|_{ L^p}$$
for any $\sigma\in\mathbb{N}$.
\et

We also introduce the bilinear version $\mathcal{B}: C^\infty(\mathbb{T}^d, \mathbb{R}^d) \times C_0^\infty(\mathbb{T}^d, \mathbb{R}^{d\times d})\to C^\infty(\mathbb{T}^d, \mathcal{S}_0^{d\times d})$ of $\mathcal{R}$.

Let $$(\mathcal{B}(v,A))_{ij}=v_l\mathcal{R}_{ijk}A_{lk}-\mathcal{R}(\partial_iv_l\mathcal{R}_{ijk}A_{lk}).$$

\bt$($\cite[Theorem B.4]{AX}$)$\label{bb}
Let $1 \leq p \leq\infty$. For any $v \in C^\infty(\mathbb{T}^d, \mathbb{R}^d)$ and $A\in C_0^\infty(\mathbb{T}^d, \mathbb{R}^{d\times d})$,

$$\div(\mathcal{B}(v,A))=vA-(2\pi)^{-d}\int_{\mathbb{T}^d}vA\dif x,$$
and
$$\|\mathcal{B}(v,A)\|_{L^p}\lesssim\|v\|_{C^1}\|\mathcal{R}A\|_{L^p}.$$
\et

\subsection{Improved H\"older’s inequality on $\mathbb{T}^d$}\label{ihiot2}
This lemma improves the usual H\"older’s inequality by using the decorrelation between frequencies.
\bt$($\cite[Theorem B.1]{AX}$)$\label{ihiot}
Let $p \in [1, \infty]$ and $a, f : \mathbb{T}^d \to \mathbb{R}$ be smooth functions. Then for any $\sigma\in\mathbb{N}$,
$$| \|a f(\sigma\cdot)\|_{L^p}-(2\pi)^{-d/p}\|a\|_{L^p}\|f\|_{L^p} |\lesssim\sigma^{-1/p}\|a\|_{C^1}\|f\|_{L^p}.$$
\et

\end{document}